\documentclass{article}
\usepackage{graphicx}

\usepackage{fullpage}
\usepackage{amsmath,amsthm,amssymb}
\usepackage{bm}
\usepackage{comment}
\usepackage{xcolor}
\usepackage{mathtools}
\usepackage{url}
\allowdisplaybreaks
\usepackage{authblk}

\usepackage{subcaption}

\title{The Inference of Fokker-Planck Equations via Transport Maps}
\author[1]{Saem Han}
\author[2]{Krishna Garikipati}
\affil[1]{Applied and Interdisciplinary Mathematics, University of Michigan}
\affil[2]{Department of Aerospace and Mechanical Engineering, University of Southern California}
\date{}

\begin{document}

\maketitle

\begin{abstract}
We present a framework, which, from the trajectories  detailing the spatiotemporal dynamics of a population,  simultaneously reconstructs a transport map as well as the Fokker-Planck equation governing the coarse-grained probability distribution. Leveraging the Knothe-Rosenblatt rearrangement, we model the transport map from a fixed reference distribution to the target distribution, and derive the velocity fields of the flows from the trajectory of transport maps. Exploiting  the velocity fields, we circumvent spatial gradients to infer the Fokker-Planck equation's potential and diffusivity. The sparsity of trajectories injects uncertainty, which we treat in a Bayesian setting using variational inference. The approach is applied to inferring the Fokker-Planck dynamics in spaces of up to five dimensions, demonstrating both accurate identification of the system and efficiency with respect to data size.
\end{abstract}

\section{Introduction}

Identifying mathematically principled representations  is centrally important to understanding the dynamics of populations. Uncovering the physical laws, and more importantly, the governing evolution equations, not only offers a rigorous understanding of the physics or mechanisms underlying the system's behavior, but also enables accurate forward predictions of the spatiotemporal dynamics. 

The literature in this broad area is rich and growing. Of relevance to our work are the studies of Brunton et al., who introduced the Sparse Identification of Nonlinear Dynamical Systems (SINDy) framework. This approach identifies ordinary or partial differential equations governing dynamical systems  from data by constructing a large dictionary of candidate functions \cite{Brunton2016}. Schaeffer proposed an optimization framework under the assumption that the underlying equations can be represented by a sparse combination of candidate terms \cite{Schaeffer2017}. Wang and co-workers developed an approach targeting the inference of the governing equations and response functions in materials systems that are dynamically evolving as well as those at mechanical equilibrium \cite{Wang2019a,wang2020system,Wang2021,wang2021system,nikolov2022ogden}. Beyond traditional frameworks that fit the equation in the strong form, they formulated the partial differential equation (PDE) identification problem in weak form. The lower requirement of smoothness on the solution field in the weak form improves robustness of inference with sparse and noisy data. Because of the weak formulation, the authors called their approach Variational System Identification. In related work, Messenger and Bortz developed a weak-form SINDy method for particle systems that demonstrated improved robustness to noise \cite{Messenger2022}. Physics-Informed Neural Networks (PINNs), introduced by Raissi et al. and further developed by Karniadakis and collaborators \cite{Raissi2019b, Karniadakis2021}, incorporate the governing PDEs directly into the loss function of neural networks, leveraging automatic differentiation to enforce physical constraints during training. This method has gained significant attention for both forward and inverse problems in scientific computing. Operator networks, in particular the DeepONet \cite{lu2019deeponet,lu2021learning} and Fourier Neural Operator approaches \cite{li2020fourier,li2023fourier} learn PDEs by constructing mappings between function spaces of input (forcing and boundary conditions) and output (solution) functions.

In this communication, we focus on the Fokker-Planck equation which describes the time evolution of a probability density function under the combined effects of random fluctuations (diffusion) and deterministic forces (drift), where the latter contribution is written as the negative gradient of a potential. The Fokker-Planck equation has played a crucial role in mathematical modeling \cite{Garrison2008,risken1989fokker}. In the context of cellular and molecular biology, Erban and Chapman demonstrated how the Fokker-Planck equation, along with stochastic simulation algorithms, can be used to model reaction-diffusion processes within cells \cite{Erban2009}. Walczak et al. developed analytic tools based on the Fokker-Planck framework to study intrinsic noise in gene expression and biochemical networks, capturing the probabilistic dynamics of single-cell behavior \cite{Walczak2012}. The Fokker-Planck equation is among a group of PDEs that have emerged as powerful tools for analysis of discretized machine learning algorithms in the continuous limit. For example, Sato and Nakagawa analyzed the stochastic gradient Langevin dynamics (SGLD) algorithm using the Fokker-Planck equation and demonstrated its convergence to the posterior average in Bayesian learning \cite{Sato2014}. Dai and Zhu investigated the statistical properties of the dynamic trajectory of stochastic gradient descent (SGD) and showed, through the Fokker-Planck framework, that the solution converges to flatter minima regardless of the batch size \cite{Dai2020}.

If the Fokker-Planck equation governs the evolution of a complex, high-dimensional, probability density function, its forward solution or inference using spatial discretization techniques, such as finite element or finite difference methods, faces limitations of exponential cost with respect to dimension. This has fueled an interest in analyzing the Fokker-Planck equation and the associated It\^{o} stochastic process in terms of simpler reference distributions and the corresponding transport maps. Going back to the seminal work of Jordan, Kinderlehrer and Otto (JKO) on the variational formulation of the Fokker-Planck equation \cite{jordan1998variational}, Liu et al. proposed a numerical method for solving Fokker-Planck equations leveraging these transport maps \cite{Liu2021}. Following JKO, these authors viewed the equation as a gradient flow of an energy functional in the quadratic Wasserstein space and introduced a family of parametric maps that represent the transport. They derived a metric tensor on the parametric space by pulling back the Wasserstein metric, and reformulated the equation as an ordinary differential equation with respect to the parameters, providing an alternative strategy generating forward solutions.

Various classes of transport maps have been studied in related settings, each offering distinct advantages. Optimal transport maps, also known as Brenier maps, minimize the quadratic cost transporting mass between two probability distributions \cite{Villani2003,Villani2008}. It is well established that the optimal transport map corresponds to the gradient of a convex function. To address the inherent computational challenges in the variational formulation, algorithms based on convex optimization have been developed. Korotin et al. proposed an approach for training optimal transport maps using the dual formulation \cite{Korotin2021a} based on the fact that the inverse mapping can be obtained as the gradient of the conjugate of the optimal potential \cite{Gangbo1996}. Makkuva et al. used convex neural networks to parametrize the convex potential whose gradient furnishes the Brenier map \cite{Makkuva2020}. 

Another important class of transport maps is the Knothe-Rosenblatt rearrangement (KR map), originally introduced in \cite{Knothe1957,Rosenblatt1952}, and further discussed in \cite{Santambrogio2015c,Villani2008}. It has gained popularity for its computational tractability - it can be computed explicitly, and the Jacobian matrix of the change of variables formula is triangular with positive entries on the diagonal. This structure has made KR map-based triangular flows a foundational component of normalizing flows in generative modeling \cite{Kobyzev2021}. Tang and Wang developed an invertible transport map by embedding the KR map into the architecture of the flow-based generative model for density estimation \cite{Tang2020}. Baptista et al. introduced an algorithm that approximates the KR map using basis functions, and proved that under certain tail conditions on the target distribution, their method converges to the unique global minimizer \cite{Baptista2023}.

In this work, our interest lies in learning the Fokker-Planck equation's driving potential and diffusion tensor. We develop a framework for an inverse problem that identifies these functions and parameters of the underlying Fokker-Planck equation governing the evolution of density data that is supplied in the form of particle trajectories at discrete time instants. Building on approaches that leverage transport maps between distributions, our method learns an explicit representation of the transport map trajectory, thereby simultaneously recovering the associated probability density functions. This step furnishes a form for the probability density flux, which also must satisfy the Fokker-Planck equation, thus presenting a convenient object to exploit for inferring the PDE. This observation is central to our approach. Section 2 reviews the Fokker-Planck equation and transport maps, and introduces the proposed methodology. In Section 3, we validate the approach using synthetic data, and make concluding remarks in Section 4.

\section{Methodology}
\begin{figure}[ht]
\includegraphics[width=0.6\textwidth]{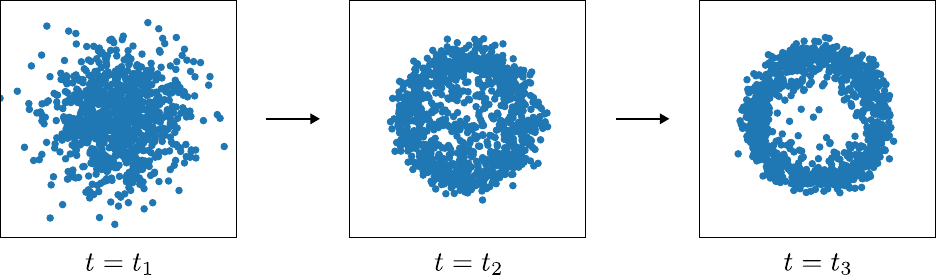}
\centering
\caption{Evolution of a data distribution governed by an It\^{o} stochastic differential equation, with snapshots  at  time instances $t=t_1, t_2, t_3$.}
\label{fig: Ito SDE}
\end{figure}

We focus on the stochastic process $X_t\in\mathbb{R}^d$ whose dynamics are governed by It\^{o} stochastic differential equations of the form: 
\begin{equation}
	dX_t = -\nabla \Psi(X_t)dt + \sigma dB_t,
    \label{eq:SDE}
\end{equation}
with the potential function $\Psi:\mathbb{R}^d\to \mathbb{R}$ and the diffusion coefficient $\sigma \in \mathbb{R}^{d\times d}$ where $B_t$ is a Brownian motion in $\mathbb{R}^d$. It is well-known that the probability density function of the distribution of the stochastic process follows the Fokker-Planck equation \cite{risken1989fokker}: 
\begin{equation}
	\frac{\partial \rho(x,t)}{\partial t} = \nabla \cdot (\nabla \Psi(x) \rho(x,t)) + \frac{1}{2}\sigma\sigma^T \nabla^2 \rho(x,t),
    \label{eq:PDE-1}
\end{equation}
where $\rho:\mathbb{R}^d\times \mathbb{R}_{\geq 0} \to \mathbb{R}$ denotes the density function of the random variable $X_t$. An example of a data distribution following this evolution equation is illustrated in Figure \ref{fig: Ito SDE}. In this paper, we aim to develop a framework for inferring the Fokker-Planck equation  that governs the dynamics of observed data, which is sample- rather than grid-based and is scalable to high-dimensional spaces. 

\subsection{Generative Learning}

The dynamics of complex distributions can be addressed by interpreting them as the result of transporting simpler reference distributions through suitable transport maps. By mapping the samples drawn from a known distribution to the target space, one can generate arbitrarily many samples to approximate the target distribution. Moreover, this framework provides an explicit representation for the target density in terms of the reference distribution and the transport map. 

We begin with two static distributions and later extend our treatment to the time-evolving setting. Let $\rho:\mathbb{R}^d\to \mathbb{R}$ and $\eta:\mathbb{R}^d\to \mathbb{R}$ denote the density functions of the target and a tractable reference distribution, respectively. Given a map $S:\mathbb{R}^d\to \mathbb{R}^d$ we say that $\rho$ is the pullback of $\eta$ under $S$ if
\begin{equation}
    \rho(x) = \eta(S(x))|\det \nabla S(x)|.
    \label{eq:pullback-1}
\end{equation}
This is expressed concisely as
\begin{equation}
    \rho = S^{\#} \eta.
    \label{eq:pullback-2}
\end{equation}
Since the transport map $S$ is not uniquely defined, learning expressive and theoretically grounded maps has been a key challenge in this approach \cite{Li2017,Arjovsky2017,Yadav2018}. One important line of research builds upon optimal transport theory \cite{Villani2003}, which adopts a variational formulation to define the map as the minimizer of the quadratic transport cost. Seguy et el. proposed a method optimizing a regularized dual form of the optimal transport cost and established its stability \cite{Seguy2018}. Subsequent works have incorporated input convex neural networks (ICNNs) \cite{Amos2017} to approximate the optimal potential functions \cite{Brenier1991,Villani2003}, leading to development of minimiax optimization algorithms \cite{Taghvaei2019,Makkuva2020}. More recently, Korotin et al. introduced a non-minimax formulation, adding a regularization term that enforces cyclic monotonicity \cite{Korotin2021a}. Recall that the quadratic Wasserstein metric is defined as: 
\begin{equation}
    W_2^2(\rho, \eta) = \min_{\rho = S^{\#} \eta} \int \frac{\|x-S(x)\|^2}{2} d\rho.
\end{equation}
According to Brenier \cite{Brenier1991} and Villani \cite{Villani2003}, the dual form is given by: 
\begin{equation}
    W_2^2(\rho, \eta) = \int \frac{\|x\|^2}{2}d\rho + \int \frac{\|y\|^2}{2}d\eta - \min_{\psi \text{: Convex}} \left[\int \psi(x)d\rho + \int \bar{\psi}(y)d\eta\right],
\end{equation}
where the minimum is taken over all  convex functions and
\begin{equation}
    \bar{\psi}(y)=\max_{x\in X} [\langle x,y\rangle - \psi(x)]
\end{equation}
is the convex conjugate to $\psi$. We obtain the following minimax formulation by parameterizing the potential function, which we denote by $\psi_\theta$: 
\begin{equation}
    \min_{\theta\in\Theta} \left[\int \psi_\theta(x)d\rho + \int \max_{x\in X} [\langle x,y \rangle - \psi_\theta(x)]d\eta\right] = \min_{\theta\in\Theta} \left[\int \psi_\theta(x)d\rho + \max_{T:Y\to X} \int [\langle T(y),y \rangle - \psi_\theta(T(y))]d\eta\right].
\end{equation}
Furthermore, we can eliminate the minimax objective by imposing additional regularization \cite{Korotin2021a}. Denote the optimal potential by $\psi^*$. As McCann \cite{McCann1995} demonstrated, using the following relation between the primal and dual potential functions:
\begin{equation}
    (\nabla \psi^*)^{-1}(y) = \nabla \bar{\psi}^*(y),
\end{equation}
we can define the regularization term as follows:
\begin{equation}
    R(\theta,w) = \int \|\nabla \psi_\theta \circ \nabla \bar{\psi}_w(y) - y\|_2^2 d\eta,
\end{equation}
where $\psi$ and $\bar{\psi}$ are parameterized using two different networks. It leads to the following non-minimax algorithm:
\begin{equation}
    \min_{\theta,w} \left[\left(\int \psi_\theta(x)d\rho + \int \left[\langle \nabla \bar{\psi}_w(y),y \rangle - \psi_\theta(\nabla \bar{\psi}_w(y))\right]d\eta\right)+\frac{\lambda}{2}R(\theta,w)\right],
\end{equation}
which imposes the condition that the optimized mappings $\nabla \psi_\theta$ and $\nabla \bar{\psi}_w$ are mutually inverse. It has been demonstrated that this non-minimax algorithm converges up to ten times faster than previous methods in \cite{Korotin2021a}. However, it may still be computationally inefficient in our framework which learns the evolution equation by iteratively updating the transport map, as the exact expression for the target density involves computing the Hessian, $\nabla^2 \psi$, which represents the Jacobian $\nabla S$, and the determinant of this Hessian as seen in the pullback expression (\ref{eq:pullback-1}).  To address this, we explore alternative approaches that can mitigate the computational cost associated with the Jacobian. 


\subsection{Knothe-Rosenblatt Rearrangements}

Among the various methods for constructing transport maps, we adopt the Knothe-Rosenblatt rearrangement \cite{Rosenblatt1952,Knothe1957}, a particularly useful choice due to its constructive definition and triangular Jacobian structure, which make it a powerful tool in mathematical modeling. Consider the two univariate densities $f(y), \, g(x):\mathbb{R}\to\mathbb{R}$ for $y\in Y$, $x \in X$, and the associated cumulative distributions $F,\,  G:\mathbb{R}\to [0,1]$.
The generalized inverse for the cumulative distribution function $F$ is defined as:
\begin{equation}
    F^{-1}(t)=\inf \left\{y \in \mathbb{R}| t \leq F(y)\right\}.
\end{equation}
A transport map between $X$ and $Y$ can be constructed in a monotone manner using the above definition of the inverse of the cumulative distribution: $S(x) = F^{-1}(G(x))$, and it is called a monotone (increasing) rearrangement. It is worth noting that this construction coincides with Brenier's optimal transportation map, which minimizes the quadratic Wasserstein distance between two measures, as discussed in \cite{McCann1995,Brenier1991,Villani2003,Villani2008}.

The components of the Knothe-Rosenblatt rearrangement are defined recursively through monotone rearrangements between the one-dimensional conditional densities. Consider the two multivariate densities $\rho, \, \eta: \mathbb{R}^d \to \mathbb{R}$. For $k = 1,\cdots,d$, let $P_k(x_k|x_1,\cdots,x_{k-1})$ and $H_k(y_k|y_1,\cdots,y_{k-1})$ denote the cumulative distribution functions of the conditional densities $\rho_k(x_k|x_1,\cdots,x_{k-1})$ and $\eta_k(y_k|y_1,\cdots,y_{k-1})$, respectively. We define the first component of $S$ as: 
\begin{equation}
    S_1(x_1) = H_1^{-1}(P_1(x_1)).
\end{equation}
Sequentially, define the $k$-th component of $S$ as: 
\begin{equation}
    S_k(x_1,\cdots,x_k) = H_k^{-1}(P_k(x_k|x_1,\cdots,x_{k-1})|S_{k-1}(x_1,\cdots,x_{k-1})).
\end{equation}
This results in a map $S:\mathbb{R}^d \to \mathbb{R}^d$, given by $S = (S_1,S_2,\cdots,S_d)^T$, which pulls back $\eta$ to $\rho$. It is triangular and monotone in the sense that each component $S_k$ depends only on the first $k$ variables, making the Jacobian matrix $\nabla S$ lower triangular, with each conditional transformation increasing in its respective variable. Consequently, $\text{det}\nabla S > 0$. The transport map satisfying these properties is uniquely defined as above and is called the Knothe-Rosenblatt rearrangement. 

Consider the following example adapted from Carlier et al. \cite{Carlier2009}, where $\eta(y)$, $y \in Y \subset \mathbb{R}^2$ is the standard Gaussian distribution with mean $\begin{pmatrix}
    0 \\ 0
\end{pmatrix}$ and covariance matrix $\begin{pmatrix}
    1 & 0 \\ 0 & 1
\end{pmatrix}$, and $\rho(x)$, $x\in X \subset \mathbb{R}^2$ is another Gaussian distribution  with 
mean $\begin{pmatrix}
    0 \\ 0
\end{pmatrix}$ and covariance matrix $\begin{pmatrix}
    a & b \\ b & c
\end{pmatrix}$, for some $a,b,c\in\mathbb{R}$. The KR map that pulls back $\eta$ to $\rho$ is:
\begin{equation}
	S(x_1,x_2) = \begin{pmatrix}
		\frac{1}{\sqrt{a}} & 0\\ -\frac{b}{\sqrt{a(ac-b^2)}}  & \frac{\sqrt{a}}{\sqrt{ac-b^2}}
	\end{pmatrix}\begin{pmatrix}
		x_1 \\ x_2
	\end{pmatrix}.
\end{equation}

It has been established (e.g., in Ref. \cite{Baptista2023}) that the KR map can be characterized as the unique minimizer of the Kullback-Leibler (KL) divergence over the space of increasing triangular maps $\mathcal{T}$:
\begin{equation}
    S = \text{arg}\,\min_{\widetilde{S} \in \mathcal{T}} D_{KL}(\rho \| \widetilde{S}^{\#}\eta)
\end{equation}
where the KL divergence quantifies the discrepancy between two distributions $\rho$ and $\widetilde{S}^{\#}\eta$ on the sample space $X$:
\begin{equation}
    D_{KL}(\rho \| \widetilde{S}^{\#}\eta) = -\sum_{x\in X} \rho(x) \log \frac{(\widetilde{S}^{\#} \eta)(x)}{\rho(x)}. 
    \label{eq:KL divergence}
\end{equation}
This approach allows us to learn the KR map that defines the pullback of the reference distribution to the target distribution, based on sample points  drawn from the target. 

\subsection{Trajectory of Transport Maps}
\label{sec: Trajectory of Transport Maps}
\begin{figure}[ht]
\includegraphics[width=0.6\textwidth]{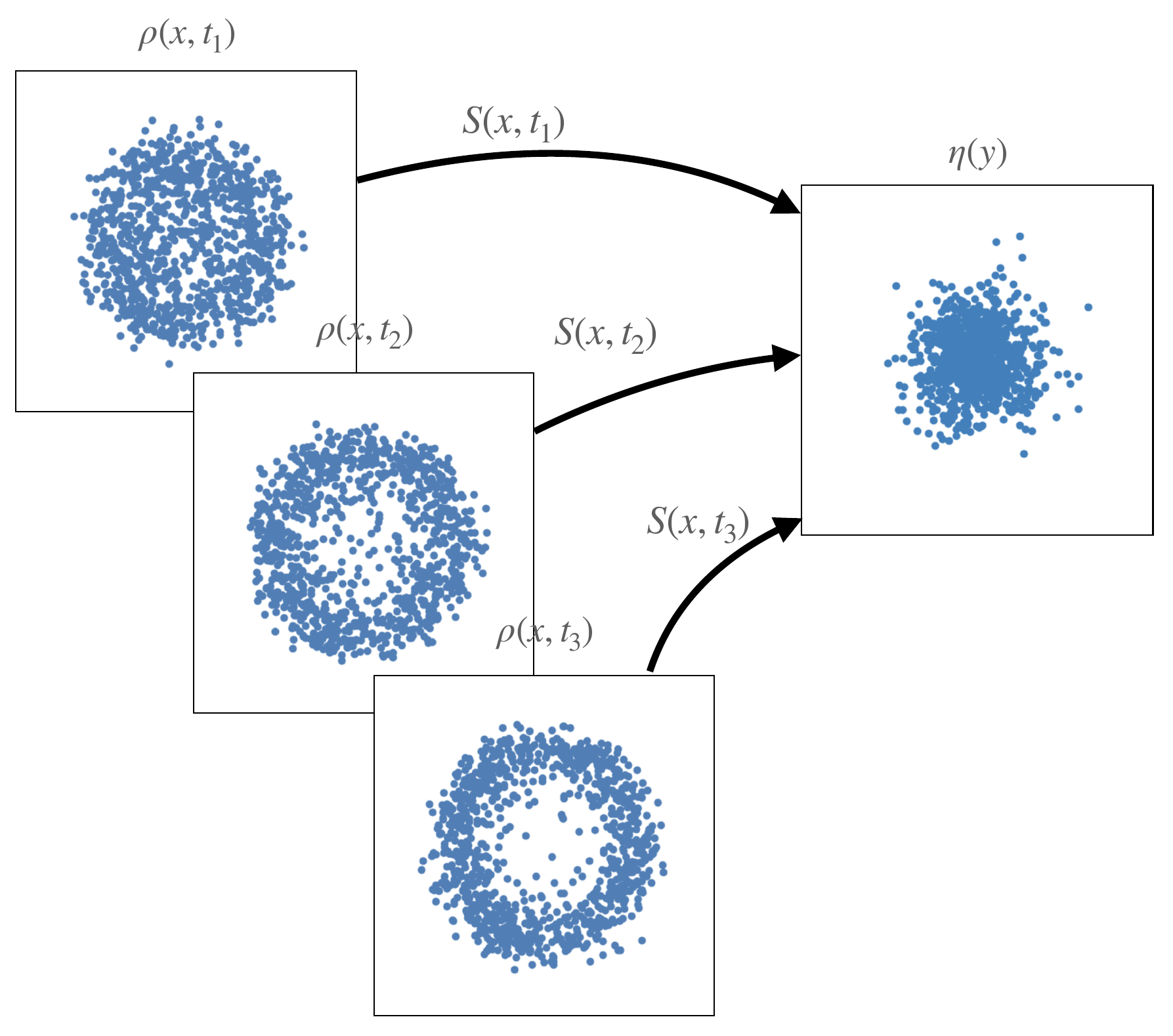}
\centering
\caption{Illustration of the transport maps $S(x,t_1)$, $S(x,t_2)$, $S(x,t_3)$ pulling back a fixed reference distribution $\eta(y)$ to target distributions $\rho(x,t_1)$, $\rho(x,t_2)$, $\rho(x,t_3)$ at different time instants. Thus, $S(x,t)$  captures the temporal evolution of the system through a sequence of diffeomorphic mappings.}
\label{fig: trajectory of mappings}
\end{figure}

We now extend our treatment to model the trajectory of maps that describe the time evolution of density functions $\rho:\mathbb{R}^d \times \mathbb{R}_{\geq 0}\to\mathbb{R}$ as pullbacks of a reference distribution $\eta: \mathbb{R}^d\to \mathbb{R}$. In this paper, we choose a standard Gaussian function for the reference distribution and assume that the target density functions are continuous and strictly positive over time. Consider the map $S:\mathbb{R}^d \times \mathbb{R}_{\geq 0} \to \mathbb{R}^d$ where $S(x,\cdot)$ is a Knothe-Rosenblatt rearrangement which pulls back $\eta(y)$ to $\rho(x,\cdot)$:
\begin{equation}
	\rho(x,t) = S(x,t)^{\#}\eta(y\circ S),
\end{equation}
where
\begin{equation}
	S(x,t) = \begin{pmatrix}
		S_1(x_1,t) \\
		S_2(x_1,x_2,t)\\
		\vdots\\
		S_d(x_1,\cdots,x_d,t)
	\end{pmatrix}.
\end{equation} 
Since the reference distribution is continuous and strictly positive, it follows that $S(x,\cdot)$ is continuously differentiable. Moreover, by the inverse function theorem, if the Jacobian has a nonzero determinant, then each KR map admits a differentiable inverse. This last condition is guaranteed since $\text{det}\nabla S$ is lower triangular with strictly positive diagonal terms. Hence, the trajectory of KR maps forms a family of diffeomorphisms. The framework is illustrated in Figure \ref{fig: trajectory of mappings}.

Since the trajectory is induced by $S(x,t)$, we can construct the velocity field of particles in the space $X$, $v:\mathbb{R}^d \times \mathbb{R}_{\geq 0} \to \mathbb{R}^d$, in a manner that characterizes the evolution relative to the static reference distribution. Using $y = S(x,t)$ and the total time derivative, yields:
\begin{equation}
	v(x,t) = -(\nabla_x S(x,t))^{-1}\frac{\partial S(x,t)}{\partial t}.
    \label{eq:velocity}
\end{equation}
We note that in the setting of continuum physics, $v(x,t)$ is called the spatial velocity field in $X$. The deformation  in continuum physics is conventionally defined as a point-to-point mapping from $Y$ to $X\times [t_0,t_N]$, including the time interval of interest in the product space. Here, it is $S^{-1}(y):\mathbb{R}^d\mapsto \mathbb{R}^d\times \mathbb{R}_{\ge 0}$.  Given the initial condition $\rho_0(x)$, the density $\rho(x,t)$ satisfies the following continuity equation with its associated initial condition:
\begin{equation}
\begin{dcases}
	\frac{\partial \rho(x,t)}{\partial t} + \nabla \cdot (\rho(x,t)v(x,t)) = 0, \quad X\times [t_0,t_N], \\
	\rho(x,t_0) = \rho_0.
\end{dcases}
\label{eq:PDE-2}
\end{equation}

Note that \eqref{eq:PDE-2} has been arrived at from conservation principles on the total mass corresponding to $\rho$. However, we also require $\rho(x,t)$ to satisfy the Fokker-Planck dynamics in \eqref{eq:PDE-1}. Our objective is to infer the underlying driving potential function and the diffusion tensor parameters that best describe the evolution of the particle distribution satisfying \eqref{eq:PDE-1} over time. Rather than directly matching two representations of the dynamics, \eqref{eq:PDE-1} and \eqref{eq:PDE-2}, we demonstrate that an equivalent result can be obtained by working solely with the corresponding fluxes. In so doing, we operate on lower-order derivatives, circumventing the divergence operator in  \eqref{eq:PDE-1}, which makes our approach inherently more resistant to noise in the data. Comparing the fluxes in \eqref{eq:PDE-1} and \eqref{eq:PDE-2}, we define an optimization problem as follows:
\begin{equation}
(\Psi,D) = 	\text{arg}\,\min_{\widetilde{\Psi}, \widetilde{D}} \left\{ \int \left| \rho v + \rho \nabla \widetilde{\Psi} +\widetilde{D}\nabla \rho \right|^2\mathrm{d}x \right\}.
    \label{eq:optimization}
\end{equation}
The associated system of Euler-Lagrange equations takes the form:
\begin{equation}
	\begin{dcases}
		\nabla \cdot(\rho (\rho v + \rho \nabla \Psi + D\nabla \rho )) = 0\\
		\nabla \rho \cdot  (\rho v + \rho \nabla \Psi + D\nabla \rho ) = 0,     
	\end{dcases}
\end{equation}
yielding
\begin{equation}
	\rho(\nabla \cdot (\rho v + \rho \nabla \Psi + D\nabla \rho))=0.
\end{equation}
Hence, under the assumption that $\rho$ is strictly positive (required for $S$ to be a diffeomorphism), the optimal $\Psi$ and $D$ satisfy:
\begin{equation}
	\nabla \cdot (\rho v) = -\nabla \cdot (\rho \nabla \Psi + D\nabla \rho). 
\end{equation}

\subsection{Model Framework}

We now model the trajectory induced by the KR maps and integrate it into the framework introduced in Section \ref{sec: Trajectory of Transport Maps}. We extend the method proposed by Baptista et al. \cite{Baptista2023}, which obtained the KR map by transforming a basis function while enforcing monotonicity. For each $k=1,\cdots, d$, let $f_k(x,t;\bm{\Theta}_1):\mathbb{R}^k \times \mathbb{R}_{\ge 0} \to \mathbb{R}$ be a smooth function with the vector of parameters $\bm{\Theta}_1$. For some positive function $g:\mathbb{R}\to\mathbb{R}_{>0}$, define the $k$-th component of the KR map as follows:
\begin{equation}
	S_k(x_1,\cdots,x_k,t; \bm{\Theta}_1) = f_k(x_1,\cdots,x_{k-1},0,t; \bm{\Theta}_1) + \int_0^{x_k} g(\partial_k f_k(x_1,\cdots,x_{k-1}, z,t; \bm{\Theta}_1))dz.
    \label{eq: Sk}
\end{equation}
This transformation ensures that $S_k$ depends only on the first $k$ spatial coordinates and is increasing with respect to the $k$-th dimension. For numerical experiments, we choose the soft-plus function $g(x) = \log(1 + e^x)$, and use standard feedforward neural networks to model each $f_k$, for their expressivity and of the framework. We collect the parameters of the Fokker-Planck potential function $\Psi$ and diffusivity $D$ in the vector $\bm{\Theta}_2$, and define the PDE-constraint as the objective function in \eqref{eq:optimization}: 
\begin{equation}
	L_{\text{PDE}}(t;\bm{\Theta}_1, \bm{\Theta}_2) = \int \left|\rho(x,t; \bm{\Theta}_1) v(x,t; \bm{\Theta}_1) + \rho(x,t; \bm{\Theta}_1)\nabla \Psi(x;\bm{\Theta}_2) + D(\bm{\Theta}_2) \nabla \rho(x,t; \bm{\Theta}_1) \right|^2 \mathrm{d}x.
    \label{eq:L_PDE}
\end{equation}
Recall that $S$ defines the target density $\rho$ via the pullback in \eqref{eq:pullback-1} and the velocity field $v$ as \eqref{eq:velocity}. For the sake of numerical stability, and to avoid the expense of repeated evaluations of $\det(\nabla S)$ and $(\nabla S)^{-1}$ we exploit \eqref{eq:pullback-1} and \eqref{eq:pullback-2} to write $\rho(x,t)v(x,t)$ as follows:
\begin{equation}
	\rho(x,t) v(x,t) = -\eta(S(x,t)) \text{adj}(\nabla_x S(x,t))\frac{\partial S(x,t)}{\partial t},
\end{equation}
where $\text{adj}(\cdot)$ is the adjugate matrix, i.e., the transpose of the cofactor matrix. 

\begin{figure}[ht]
\includegraphics[width=0.6\textwidth]{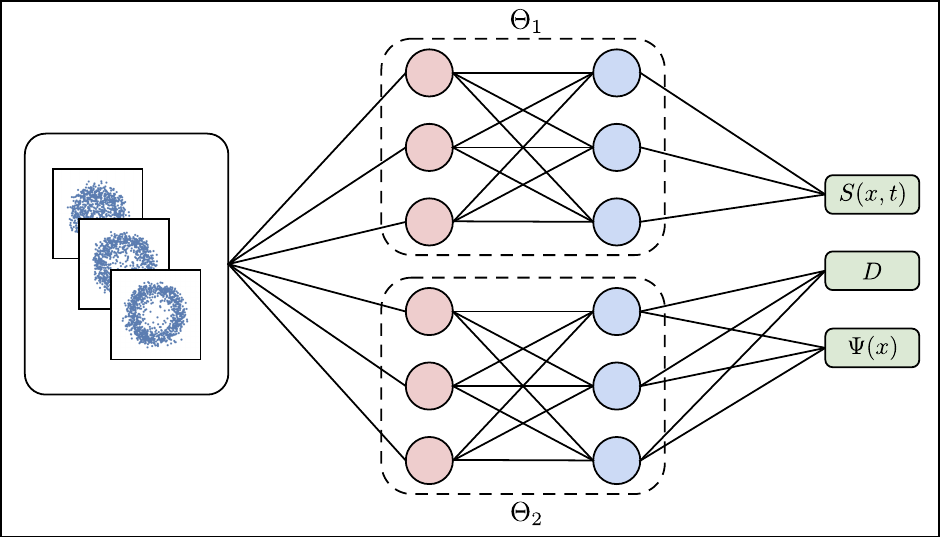}
\centering
\caption{Overview of the model architecture. The pullback map $S$ and the Fokker-Planck equation quantities $\Psi,D$ are parametrized by $\bm{\Theta}_1$ and $\bm{\Theta}_2$, respectively, and optimized based on data observed at discrete time points.}
\label{fig: model architecture}
\end{figure}

Given data at discrete time stamps $t=t_i$, $i=0,\cdots,N$, we define the loss function by combining the KL divergence \eqref{eq:KL divergence} with the PDE-constraint \eqref{eq:L_PDE}, weighted by a tunable parameter $\lambda$:
\begin{equation}
	L = \sum_{i=1}^N \left\{ D_{KL} (\rho(x,t_i; \bm{\Theta}_1)\|S(x,t_i; \bm{\Theta}_1)^\# \eta(x)) + \lambda L_{\text{PDE}}(t_i;\bm{\Theta}_1, \bm{\Theta}_2) \right\}.
    \label{eq: deterministic loss - KL}
\end{equation}
Given a sample space $X$, the KL divergence in \eqref{eq: deterministic loss - KL} decomposes as: 
\begin{equation}
    D_{KL}(\rho \| S^{\#}\eta) = -\sum_{x\in X} \left(\rho(x) \log (S^{\#} \eta)(x) - \rho(x)\log \rho(x)\right). 
    \label{eq:KL divergence 2}
\end{equation}
The second term in \eqref{eq:KL divergence 2} is entirely data-dependent, and is a constant for the optimization problem. We rewrite the objective function \eqref{eq: deterministic loss - KL} by replacing the KL divergence with the negative log-likelihood as follows:
\begin{equation}
	L = \sum_{i=1}^N \left\{ -\mathbb{E}_{\rho} [\log S(x,t_i;\bm{\Theta}_1)^{\#}\eta] + \lambda L_{\text{PDE}}(t_i;\bm{\Theta}_1, \bm{\Theta}_2) \right\},
    \label{eq: deterministic loss}
\end{equation}
where $\mathbb{E}_\rho$ is the expectation over $\rho$. The log-likelihood function is evaluated over the given sample space
\begin{equation}
    \mathbb{E}_{\rho} [\log S^{\#}(x,t_i;\bm{\Theta}_1)\eta] = \sum_{x\sim \rho(x,t_i)} \log S^{\#}(x,t_i;\bm{\Theta}_1) \eta(x),
\end{equation}
and $L_{\text{PDE}}$ can be approximated using Monte-Carlo methods with random sampling. By training the model with the loss function \eqref{eq: deterministic loss}, we simultaneously learn the continuous representation of the map $S(x,t;\bm{\Theta}_1)$ and the system parameters $\bm{\Theta}_2$. The model architecture is illustrated in Figure \ref{fig: model architecture}.

\subsubsection{Uncertainty quantification}
The finite data introduces uncertainty in our inference of the transport map and Fokker-Planck functions and parameters. We extend the framework to a probabilistic setting that interprets each parameter as a sample from an underlying distribution. Suppose that the set of parameters $\tilde{\bm{\Theta}}_1$ and $\tilde{\bm{\Theta}}_2$ are used to describe such distributions for the Knothe-Rosenblatt rearrangements and the Fokker-Planck equation, respectively.

By Bayes' theorem, $P(\tilde{\bm{\Theta}}_j|\mathcal{D})$, the posterior distribution of the parameters given data $\mathcal{D}$, is written as:
\begin{equation}
    P(\tilde{\bm{\Theta}}_j|\mathcal{D}) = \frac{P(\mathcal{D}|\tilde{\bm{\Theta}}_j)P(\tilde{\bm{\Theta}}_j)}{P(\mathcal{D})},
\end{equation}
where $P(\mathcal{D}|\tilde{\bm{\Theta}}_j)$  and $P(\tilde{\bm{\Theta}}_j)$ denote the likelihood and the prior probability distributions, respectively, for $j=1,2$, and $P(\mathcal{D})$ represents the evidence. Since computing $P(\mathcal{D})$ is intractable as it requires marginalization over all parameters, we approximate the posterior distribution with a more tractable surrogate distribution, denoted by $Q_j$, by minimizing the KL divergence
$D_{KL}(Q_j(\tilde{\bm{\Theta}}_j)\| P(\tilde{\bm{\Theta}}_j|\mathcal{D}))$. We follow the steps adopted widely in the literature \cite{Blei2017,Graves2011,Zhang2023}. We observe also that the KL divergence is equal to the negative evidence lower bound (ELBO) up to the constant evidence term
\begin{align}
    D_{KL}(Q_j(\tilde{\bm{\Theta}}_j)\| P(\tilde{\bm{\Theta}}_j|\mathcal{D})) &= - \underbrace{\left\{\mathbb{E}_{Q_j}[\log P(\tilde{\bm{\Theta}}_j, \mathcal{D})] - \mathbb{E}_{Q_j}[\log Q_j (\tilde{\bm{\Theta}}_j)]  \right\} }_{ELBO} + \log P(\mathcal{D})\\
    &=D_{KL}(Q_j(\tilde{\bm{\Theta}}_j)\| P(\tilde{\bm{\Theta}}_j)) - \mathbb{E}_{Q_j}[\log P(\mathcal{D}|\tilde{\bm{\Theta}}_j)] + \log P(\mathcal{D}).
\end{align}
Therefore, we define a probabilistic loss as follows:
\begin{equation}
	L_\text{prob} = \left\{ D_{KL}(Q_1(\tilde{\bm{\Theta}}_1)\|P(\tilde{\bm{\Theta}}_1)) - \mathbb{E}_{Q_1}[\log P(\mathcal{D}|\tilde{\bm{\Theta}}_1)] \right\} + \lambda \left\{ D_{KL}(Q_2(\tilde{\bm{\Theta}}_2)\|P(\tilde{\bm{\Theta}}_2)) - \mathbb{E}_{Q_2}[\log P(\mathcal{D}|\tilde{\bm{\Theta}}_2)] \right\}.
    \label{eq: probablistic loss}
\end{equation}
 As with the deterministic loss \eqref{eq: deterministic loss} the first two terms parameterized by $\tilde{\bm{\Theta}}_1$  learn the transport map, while the last two parameterized by $\tilde{\bm{\Theta}}_2$  enforce the Fokker-Planck PDE. The expectations of log-likelihoods, by definition, control the representation of the data by way of samples, and by the PDE, respectively. The KL-divergences provide regularization by keeping the surrogate posteriors close to the priors.
 
 In practice, the deterministic model can be pre-trained to provide an effective initialization for the Bayesian model, as demonstrated by Zhang and Garikipati \cite{Zhang2023}. Furthermore, the learned deterministic parameters $\bm{\Theta}_1$ can be re-used in the log-likelihood term within the PDE constraint, expressing $P(\mathcal{D}|\tilde{\bm{\Theta}}_2)$ as $P(\mathcal{D};\bm{\Theta}_1|\tilde{\bm{\Theta}}_2)$, to help stabilize Bayesian training:
 \begin{equation}
	L_\text{prob} = \left\{ D_{KL}(Q_1(\tilde{\bm{\Theta}}_1)\|P(\tilde{\bm{\Theta}}_1)) - \mathbb{E}_{Q_1}[\log P(\mathcal{D}|\tilde{\bm{\Theta}}_1)] \right\} + \lambda \left\{ D_{KL}(Q_2(\tilde{\bm{\Theta}}_2)\|P(\tilde{\bm{\Theta}}_2)) - \mathbb{E}_{Q_2}[\log P(\mathcal{D};\bm{\Theta}_1|\tilde{\bm{\Theta}}_2)] \right\}.
    \label{eq: probablistic loss}
\end{equation}

The likelihood functions can be specified heuristically and parameterized accordingly. Since $P(\mathcal{D}|\tilde{\bm{\Theta}}_1)$ assigns a higher probability when the KL divergence is close to zero, it can be modeled as a Gaussian distribution with mean $0$ and trainable variance $\sigma_1$: $P(\mathcal{D}|\tilde{\bm{\Theta}}_1) = \mathcal{N}\left(\sum_{i} D_{KL} (\rho(x,t_i; \tilde{\bm{\Theta}}_1)\|S^\#(x,t_i; \tilde{\bm{\Theta}}_1) \eta(x))  \bigg|0, \sigma_1\right)$. Another likelihood $P(\mathcal{D}; \bm{\Theta}_1|\tilde{\bm{\Theta}}_2)$ enforces the PDE-constraint which minimizes \eqref{eq:L_PDE}. Denoting the minimum of \eqref{eq:L_PDE} as $L^*$, it can be represented as a Gaussian distribution with trainable variance $\sigma_2$: $P(\mathcal{D}; \bm{\Theta}_1|\tilde{\bm{\Theta}}_2) = \mathcal{N}\left(\sum_{i}  L_\text{PDE}(t_i;\bm{\Theta}_1, \tilde{\bm{\Theta}}_2)\bigg|L^*, \sigma_2\right)$. Incorporating these steps, we have:

\begin{align}
    L_{\text{prob}} &= \left\{ D_{KL}(Q_1(\tilde{\bm{\Theta}}_1)\|P(\tilde{\bm{\Theta}}_1)) - \mathbb{E}_{Q_1}\left[\log \mathcal{N}\left(\sum_{i} D_{KL} (\rho(x,t_i; \tilde{\bm{\Theta}}_1)\|S^\#(x,t_i; \tilde{\bm{\Theta}}_1) \eta(x))  \bigg|0, \sigma_1\right)\right] \right\} \nonumber\\
    &+ \lambda \left\{ D_{KL}(Q_2(\tilde{\bm{\Theta}}_2)\|P(\tilde{\bm{\Theta}}_2)) - \mathbb{E}_{Q_2}\left[\log \mathcal{N}\left(\sum_{i}  L_\text{PDE}(t_i;\bm{\Theta}_1, \tilde{\bm{\Theta}}_2)\bigg|L^*, \sigma_2\right)\right] \right\}.
    \label{eq:finprobloss}
\end{align}

\section{Computational Results}
We train and test our inference framework against synthetic data generated by numerically solving stochastic differential equations \eqref{eq:SDE} using the Euler-Maruyama method with a sufficiently small time step to ensure high accuracy \cite{Talay1994}. The particle distributions are obtained at discrete time instances $t_i = 0.1i$, for $i=0,1,\cdots, 10$, to simulate the temporal evolution of the system. We generate a varying number of sample points depending on the dimensionality, allowing more samples in higher-dimensional settings to ensure adequate representation of the underlying distribution. We show that the model successfully learns the map $S$, accurately reconstructs the density function, and correctly identifies the underlying system parameters.

The core ingredient of the pullback map \eqref{eq: Sk}, $f_k$, is represented by a feed-forward neural network with five hidden layers, each containing five nodes, and utilizing the hyperbolic tangent activation function. Although the hyperparameters were determined by a grid search on a set of sample problems and kept fixed across all experiments presented in this paper, we anticipate that problem-specific tuning could yield more optimal results. For the reported experiments, the following hyperparameter settings were used: The $\lambda$ in the loss functions \eqref{eq: deterministic loss} and \eqref{eq:finprobloss} was set to $0.1$. The model was trained using the Adam optimizer with a learning rate $0.001$ and decay rate parameters $\beta_1 = 0.8$ and $\beta_2=0.999$. The partial derivative with respect to the last spatial variable, $\partial_k f_k$, is calculated exactly and efficiently using the built-in backpropagation functions in \texttt{PyTorch} for $k = 1,\dots, d$. The corresponding one-dimensional definite integral is then approximated numerically using Simpson's rule with $20$ partitions along the axis.

\subsection{Two-dimensional Problems}
We consider the initial distribution $\rho_0$ to be a mixture of four Gaussian components:
\begin{equation}
    \rho_0(x) = \frac{1}{4}\sum_{i=1}^4 \mathcal{N}(x\vert\mu_i, \Sigma_i),
    \label{eq: Initial Condition - MG}
\end{equation}
where the means and covariances of the components are specified to be:
\begin{equation}
    \mu_1 = \begin{pmatrix}
2 \\ 2 \end{pmatrix}, \;
\mu_2 = \begin{pmatrix}
-2 \\ 2 \end{pmatrix}, \;
\mu_3 = \begin{pmatrix}
-2 \\ -2 \end{pmatrix},\; 
\mu_4 = \begin{pmatrix}
2 \\ -2 \end{pmatrix}, 
\end{equation}
and 
\begin{equation}
    \Sigma_1 = \begin{pmatrix}
    1 & 0 \\ 0 & 1
\end{pmatrix},\;
\Sigma_2 = \begin{pmatrix}
    2 & 0 \\ 0 & 0.5
\end{pmatrix},\;
\Sigma_3 = \begin{pmatrix}
    0.5 & 0 \\ 0 & 1
\end{pmatrix},\;
\Sigma_4 = \begin{pmatrix}
    1 & 0 \\ 0 & 2
\end{pmatrix}.
\end{equation}

\subsubsection{Isotropic diffusion with quadratic potential functions}

\begin{figure}[ht!]
\includegraphics[width=0.8\textwidth]{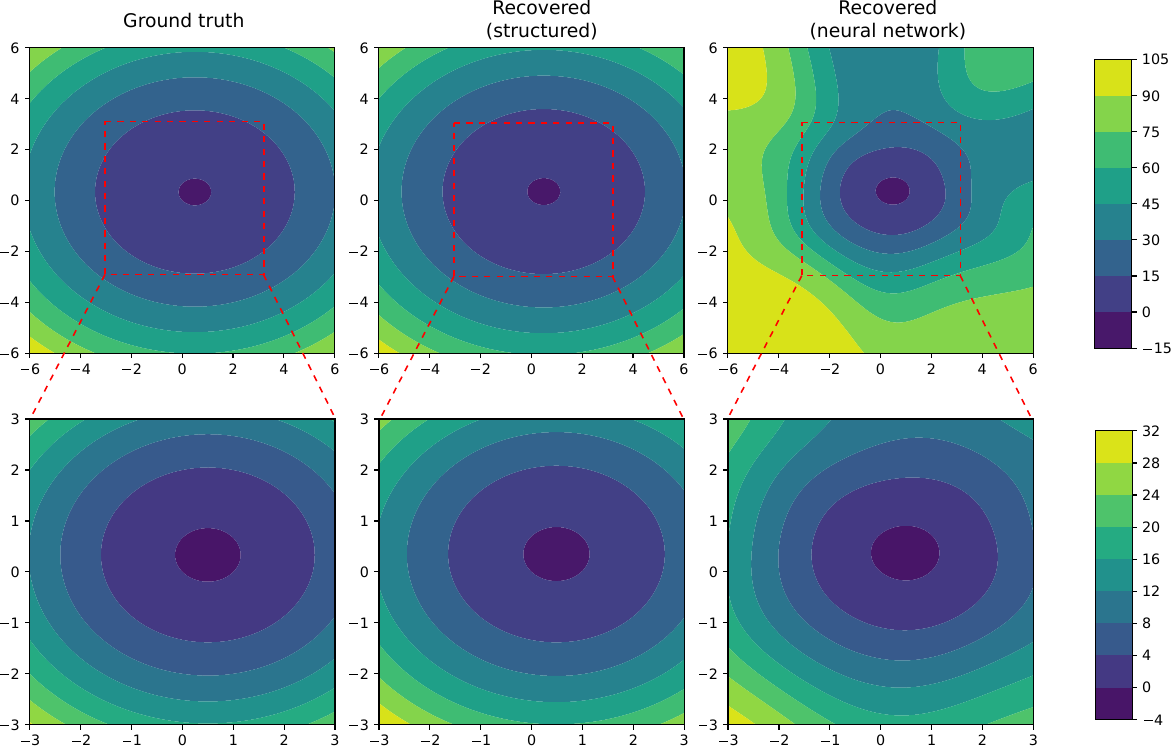}
\centering
\caption{Contour plots of the potential function for the two-dimensional isotropic, quadratic potential problem. The lower row is the zoomed-in version of the upper row images (see axis limits). Each column shows: (left) the ground truth potential, (middle) reconstructed potential using the structured model, and (right) the reconstructed potential using the neural network model. }
\label{fig: 2D Isotropic Quadratic - Potential}
\end{figure}

\begin{figure}[ht!]
\includegraphics[width=0.7\textwidth]{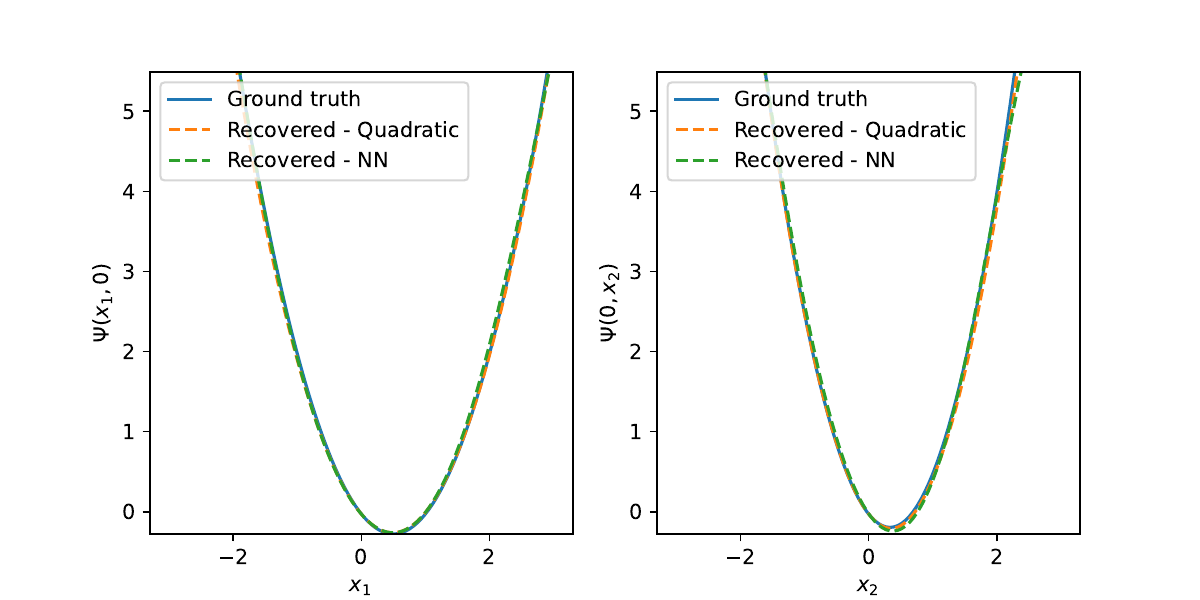}
\centering
\caption{Slices of the potential function along each axis: (left) at $x_2=0$ and (right) at $x_1=0$ for the two-dimensional isotropic, quadratic potential problem.}
\label{fig: 2D Isotropic Quadratic - Potential All}
\end{figure}

\begin{figure}[ht!]
    \centering
    \includegraphics[width=0.6\textwidth]{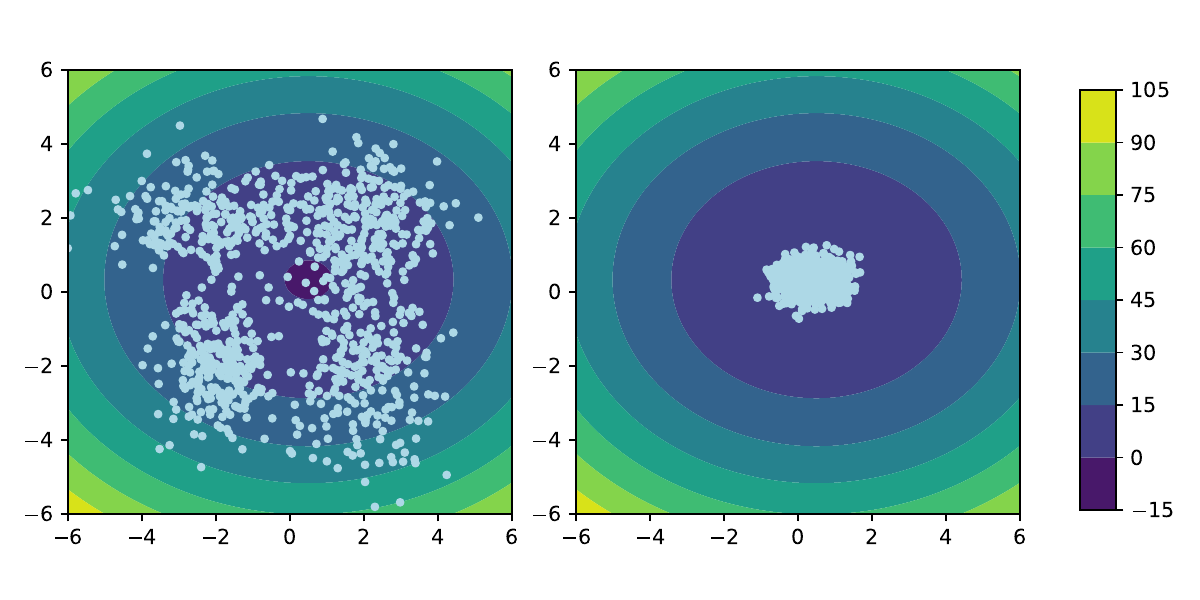}
    \caption{The initial (left) and the final (right) data samples (light blue dots) superposed on the contour plot of the Fokker-Planck potential function for the two-dimensional isotropic, quadratic potential problem.}
    \label{fig: 2D Isotropic Quadratic - Potential and Samples}
\end{figure}

Suppose that the  density evolves from the initial distribution $\rho_0$  according to the
Fokker-Planck equation with an isotropic diffusion tensor and a quadratic potential functions shown below:
\begin{equation}
	D =\begin{pmatrix}
 0.2 & 0 \\ 0 & 0.2	
 \end{pmatrix}, \quad
 \Psi(x) = \frac{1}{2}x^T\begin{pmatrix}
 	2 & 0 \\ 0 & 3
 \end{pmatrix}x - \begin{pmatrix}
 	1 \\ 1
 \end{pmatrix}^T x. 
\end{equation}
For this experiment, $1,000$ sample points are drawn from the initial distribution, and the entire dataset is constructed by numerically solving the associated stochastic differential equation by the Euler-Maruyama method. The diffusion tensor is modeled via a scalar-valued parameter, reflecting uniform diffusion in all directions:
\begin{equation}
    D =\begin{pmatrix}
 \theta_D & 0 \\ 0 & \theta_D	
 \end{pmatrix}.
\end{equation}
For the potential function, we explore two approaches. In the first, we assume prior knowledge about the potential function's structure and use a quadratic form:
\begin{equation}
    \Psi(x) = \frac{1}{2}x^T\begin{pmatrix}
 	\theta_\Psi^1 & 0 \\ 0 & \theta_\Psi^2
 \end{pmatrix}x + \begin{pmatrix}
 	\theta_\Psi^3 \\ \theta_\Psi^4
 \end{pmatrix}^T x.
\end{equation}
In the second approach, we relax this assumption and represent the potential function as a learnable neural network. The potential functions reconstructed using these two  models appear in Figure \ref{fig: 2D Isotropic Quadratic - Potential} and Figure \ref{fig: 2D Isotropic Quadratic - Potential All}. Over a wider domain, as shown in the top row of Figure \ref{fig: 2D Isotropic Quadratic - Potential}, the neural network potential model, not being constrained to the ground-truth quadratic form, exhibits noticeably different behavior, which is expected since the data are concentrated closer to the origin (see Figure \ref{fig: 2D Isotropic Quadratic - Potential and Samples}). Within the region where the data points are clustered, depicted in the bottom row of Figure \ref{fig: 2D Isotropic Quadratic - Potential}, both models accurately capture the underlying potential function.  The diffusion tensor and potential parameters obtained under the assumption of structured form are summarized in Table \ref{table: 2D Isotropic Quadratic}.

\begin{table}[ht!]
       \centering
       \begin{tabular}{|c|c|c|c|c|c|}
       \hline 
            & $\theta_D$ & $\theta_\Psi^1$ & $\theta_\Psi^2$ & $\theta_\Psi^3$ & $\theta_\Psi^4$\\
            \hline
            Ground truth & 0.2 & 2.0 & 3.0 & -1.0 & -1.0\\
           \hline
           Recovered & 0.19 & 1.97 & 2.94 & -0.97 & -1.03\\
           \hline
       \end{tabular}
       \caption{Ground truth and recovered system parameters for the two-dimensional isotropic, quadratic potential problem. }
       \label{table: 2D Isotropic Quadratic}
 \end{table}
 
\begin{figure}[ht!]
\includegraphics[width=1.0\textwidth]{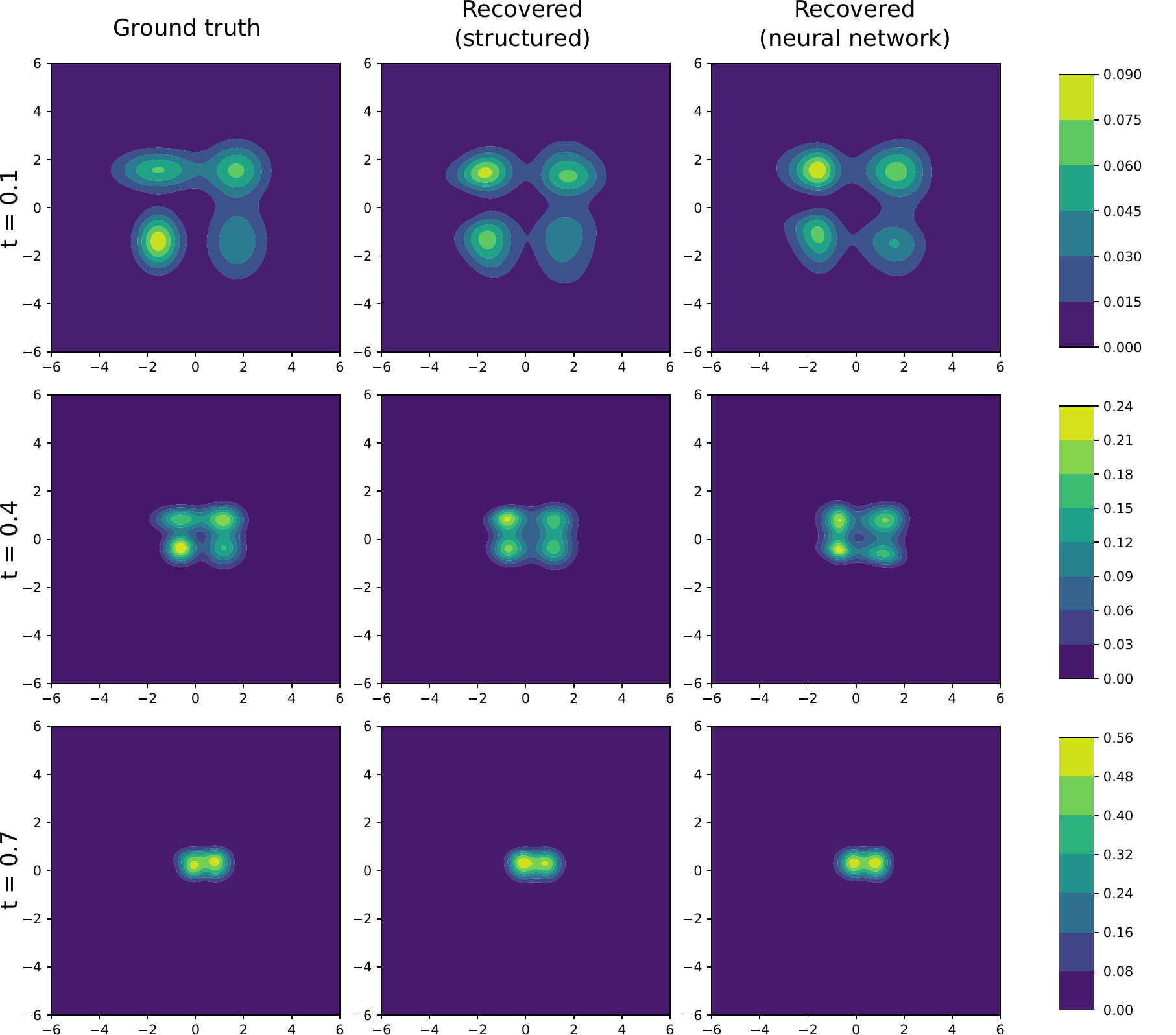}
\centering
\caption{Contour plots of the density at different time instances. Each row corresponds to a different time stamp: $t=0.1, 0.4, 0.7$ (from top to bottom). Each column shows: (left) ground truth density, (middle) reconstructed density using the structured potential model, and (right) reconstructed density using the neural network potential model.}
\label{fig: 2D Isotropic Quadratic - PDF}
\end{figure}

Furthermore, the continuous representation of the KR map trajectory is obtained from training, allowing the density function to be explicitly evaluated by \eqref{eq:pullback-1}. Figure \ref{fig: 2D Isotropic Quadratic - PDF} shows the reconstructed probability densities at selected time points after training. The ground truth density functions, illustrated in the leftmost column, are generated from the analytic solution of the inferred Fokker-Planck equation. 

\subsubsection{Anisotropic diffusion with annular potential functions}

\begin{figure}[ht!]
\includegraphics[width=0.8\textwidth]{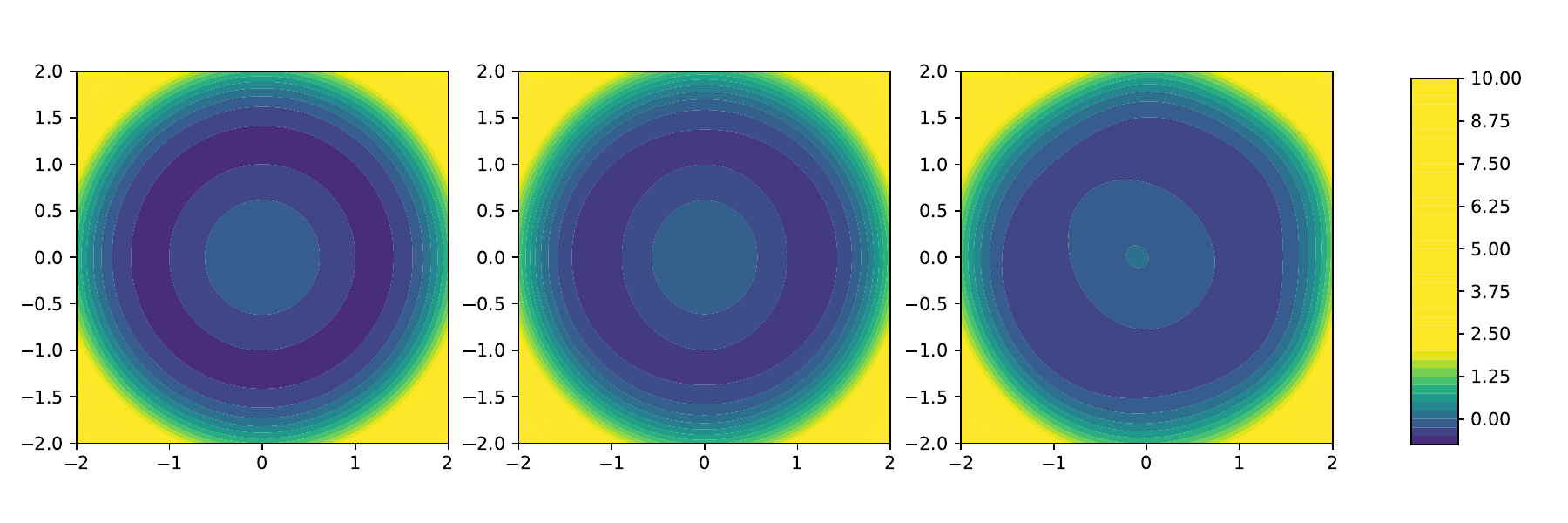}
\centering
\caption{Contour plots of the potential function. Each column shows: (left) the ground truth potential, (middle) reconstructed potential using the structured model, and (right) the reconstructed potential using the neural network model for the two-dimensional anisotropic, annular potential problem.}
\label{fig: 2D Anisotropic DoubleWell - Potential}
\end{figure}

\begin{figure}[ht!]
\includegraphics[width=0.7\textwidth]{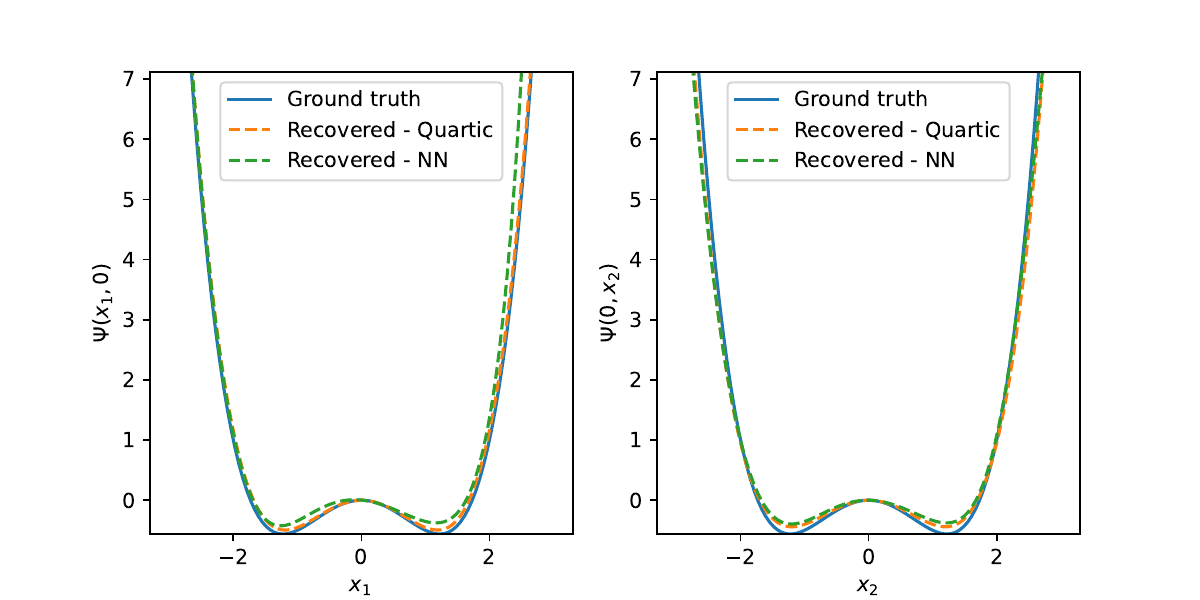}
\centering
\caption{Slices of the potential function along each axis: (left) at $x_2=0$ and (right) at $x_1=0$ for the two-dimensional anisotropic, annular potential problem. }
\label{fig: 2D Anisotropic DoubleWell - Potential Slices}
\end{figure}

\begin{figure}[ht!]
    \centering
    \includegraphics[width=0.6\linewidth]{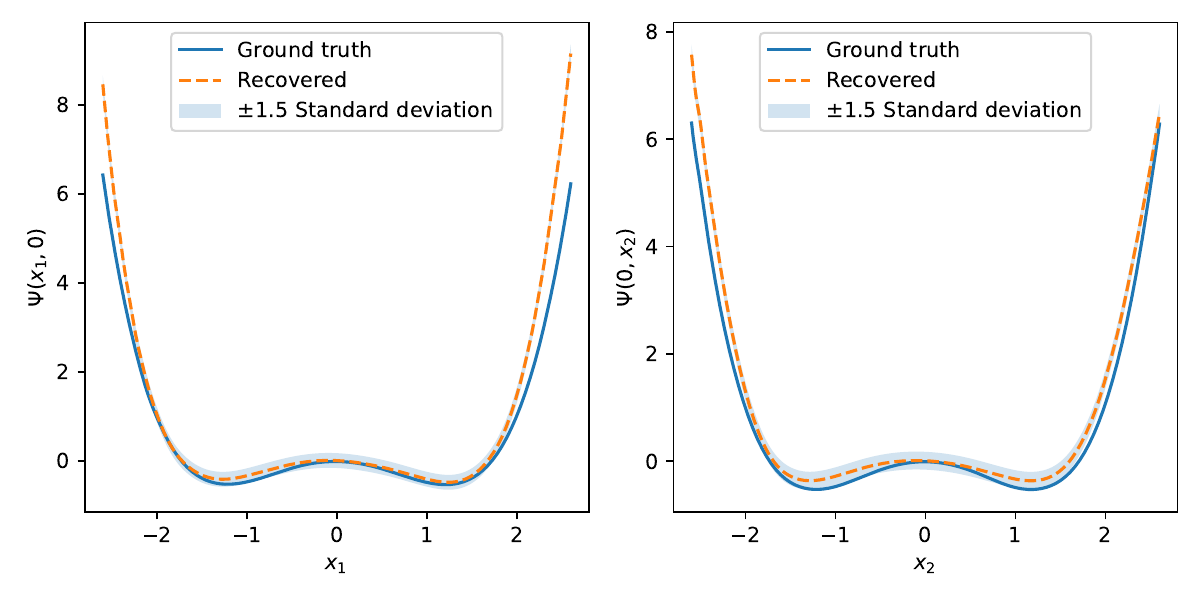}
    \caption{Slices of the potential function reconstructed by the Bayesian framework. The profiles are shown along each axis: (left) at $x_2=0$ and (right) at $x_1=0$ for the two-dimensional anisotropic, annular potential problem. The dashed lines represent the mean, while the shaded regions illustrate the standard deviation.}
    \label{fig: 2D Anisotropic DoubleWell Potential - BNN}
\end{figure}

\begin{figure}[ht!]
    \centering
    \includegraphics[width=0.6\textwidth]{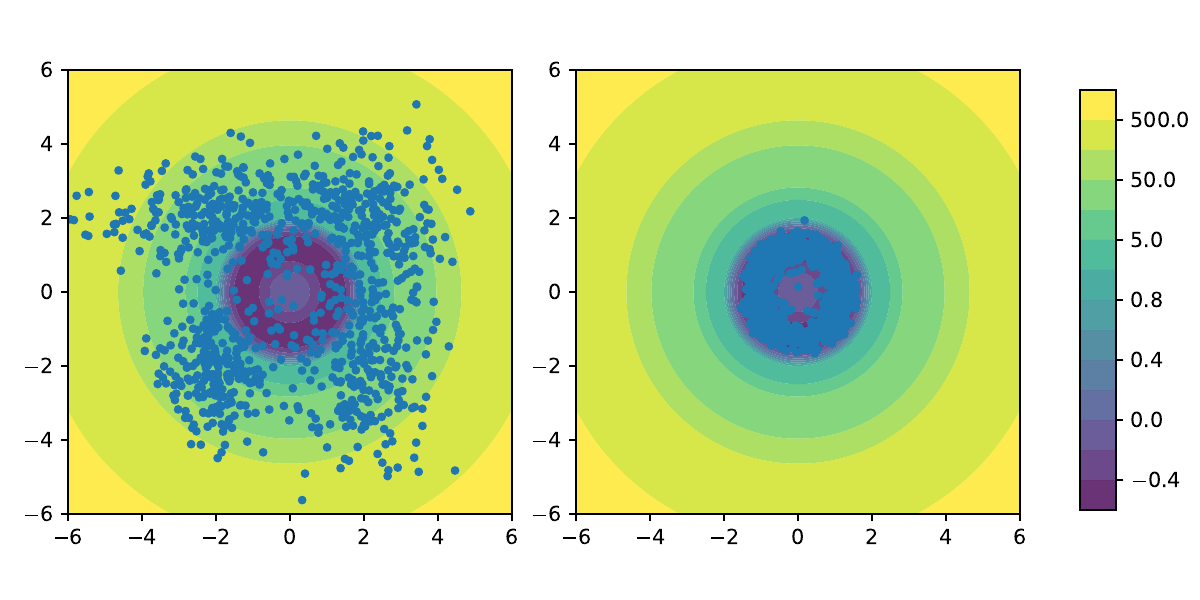}
    \caption{The initial (left) and the final (right) data samples (blue dots) superposed on the contour plot of the Fokker-Planck potential function for the two-dimensional anisotropic, annular potential problem. }
    \label{fig: 2D Anisotropic DoubleWell - Potential and Samples}
\end{figure}

We next consider dynamics which are governed by anisotropic diffusion and the annular potential function:
\begin{equation}
	D =\begin{pmatrix}
 0.1 & 0 \\ 0 & 0.2	
 \end{pmatrix}, \quad
 \Psi(x) = \frac{1}{4}\left(x^T\begin{pmatrix}
 	1 & 0 \\ 0 & 1
 \end{pmatrix}x\right)^2 -\frac{1}{2}x^T\begin{pmatrix}
 	1.5 & 0 \\ 0 & 1.5
 \end{pmatrix}x 
\end{equation}
sampling from a mixture of four Gaussians for $\rho_0$ in \eqref{eq: Initial Condition - MG} as the initial condition. Again, the training set is formed by using $1,000$ data points drawn from the initial distribution. For anisotropic diffusion problems, we allow a full matrix-valued structure for the diffusion model, enabling directionally dependent diffusion behaviors to be learned from the data:
\begin{equation}
    D =\begin{pmatrix}
 \theta_D^1 & \theta_D^2 \\ \theta_D^2 & \theta_D^3	
 \end{pmatrix}
\end{equation}
for symmetric diffusion tensors. The potential function is then modeled in two ways as in the previous problem, first, with the pre-defined structure:
\begin{equation}
    \Psi(x) = \frac{1}{4}\left(x^T\begin{pmatrix}
 	\theta_\Psi^1 & 0 \\ 0 & \theta_\Psi^2
 \end{pmatrix}x\right)^2 + \frac{1}{2}x^T\begin{pmatrix}
 	\theta_\Psi^3 & 0 \\ 0 & \theta_\Psi^4
 \end{pmatrix}x,
\end{equation}
and second, using a general neural network. The recovered potential functions along with the ground truth are shown in Figure \ref{fig: 2D Anisotropic DoubleWell - Potential} and Figure \ref{fig: 2D Anisotropic DoubleWell - Potential Slices}. The structured model clearly recovers the underlying potential function more accurately, while the neural network model, despite being unconstrained, successfully captures the shape of the wells. Figure \ref{fig: 2D Anisotropic DoubleWell Potential - BNN} shows the mean and standard deviation of the Bayesian neural network-learned potential plotted in comparison with the ground truth. Table \ref{table: 2D Anisotropic DoubleWell} presents the inferred diffusion and potential parameters under the structure-form assumption. The initial and final data samples superposed on the annular potential function are displayed in Figure \ref{fig: 2D Anisotropic DoubleWell - Potential and Samples}.  

\begin{table}[ht!]
       \centering
       \begin{tabular}{|c|c|c|c|c|c|c|c|c|c|c|c|c|}
       \hline 
            & $\theta_D^1$ & $\theta_D^2$  & $\theta_D^3$ & $\theta_\Psi^1$ & $\theta_\Psi^2$ & $\theta_\Psi^3$ & $\theta_\Psi^4$ \\
            \hline
            Ground truth & 0.1 & 0.0 & 0.2 & 1.0 & 1.0 & -1.5 & -1.5 \\
           \hline
           Recovered & 0.11 & -0.01 & 0.17 & 1.01 & 0.91 & -1.45 & -1.21 \\
           \hline
       \end{tabular}
       \caption{Ground truth and recovered system parameters for the two-dimensional anisotropic, annular potential problem.}
       \label{table: 2D Anisotropic DoubleWell}

\end{table}

\begin{figure}[ht!]
\includegraphics[width=1.0\textwidth]{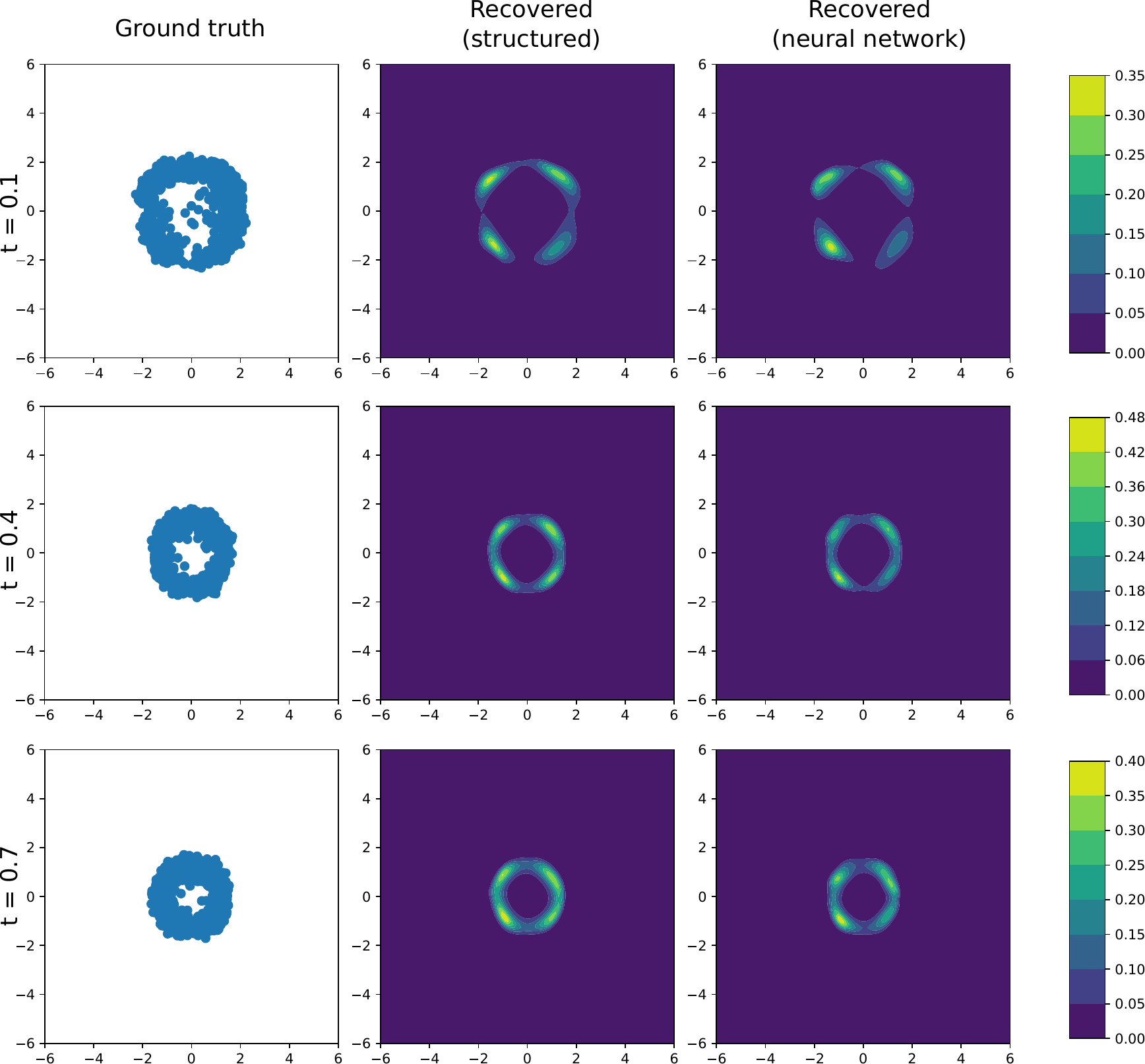}
\centering
\caption{Contour plots of the density at different time instances. The rows correspond to a different time stamp: $t=0.1, 0.4, 0.7$ (from top to bottom). Each column shows: (left) ground truth density, (middle) reconstructed density using the structured potential model, and (right) reconstructed density using the neural network potential model for the two-dimensional anisotropic, annular potential problem.}
\label{fig: 2D Anisotropic DoubleWell - PDF}
\end{figure}

\begin{figure}[ht!]
\includegraphics[width=1.0\textwidth]{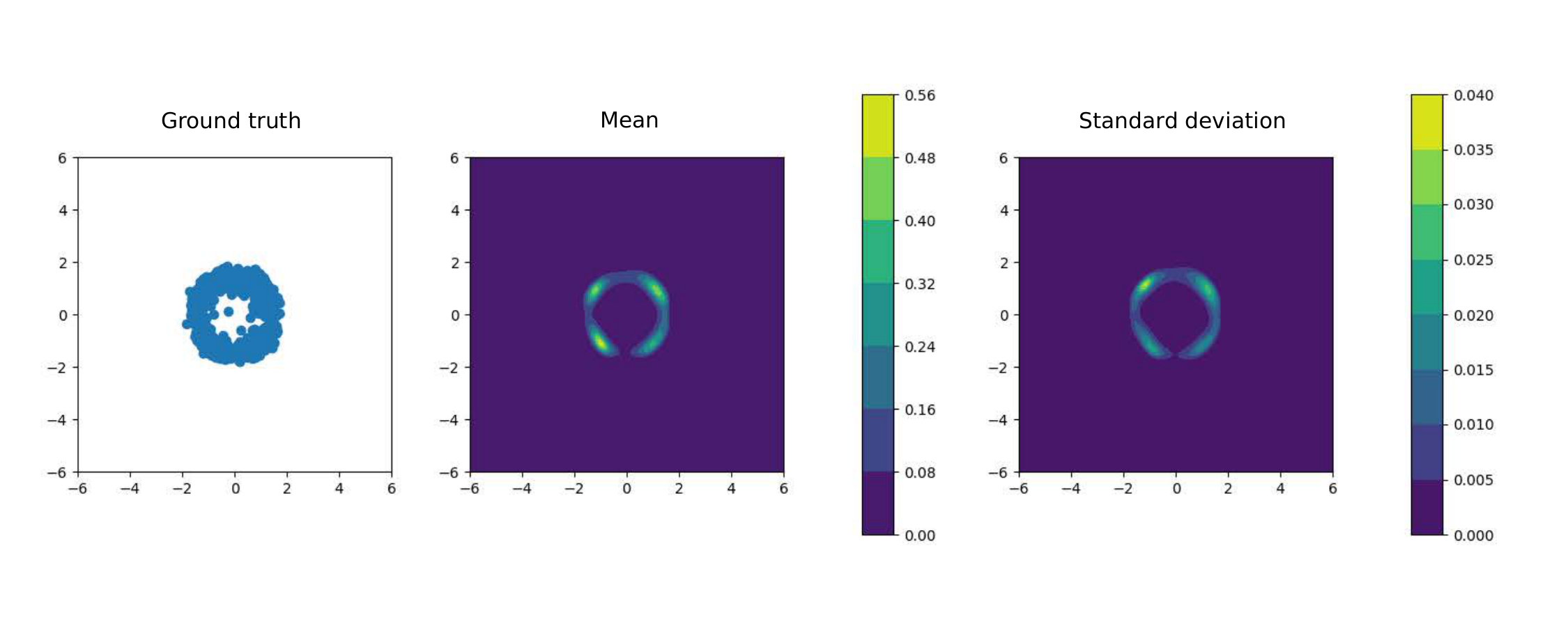}
    \centering
    \caption{Mean and standard deviation of the reconstructed density at $t=1.0$ generated by bayesian neural networks with a probabilistic architecture. (Left) ground truth distribution, (middle) mean of the reconstructed density, and (right) standard deviation of the reconstructed density obtained using the neural network potential model for the two-dimensional anisotropic, annular potential problem.}
    \label{fig: 2D Anisotropic DoubleWell - BNN}
\end{figure}

\begin{figure}[ht!]
    \centering
    \includegraphics[width=0.36\linewidth]{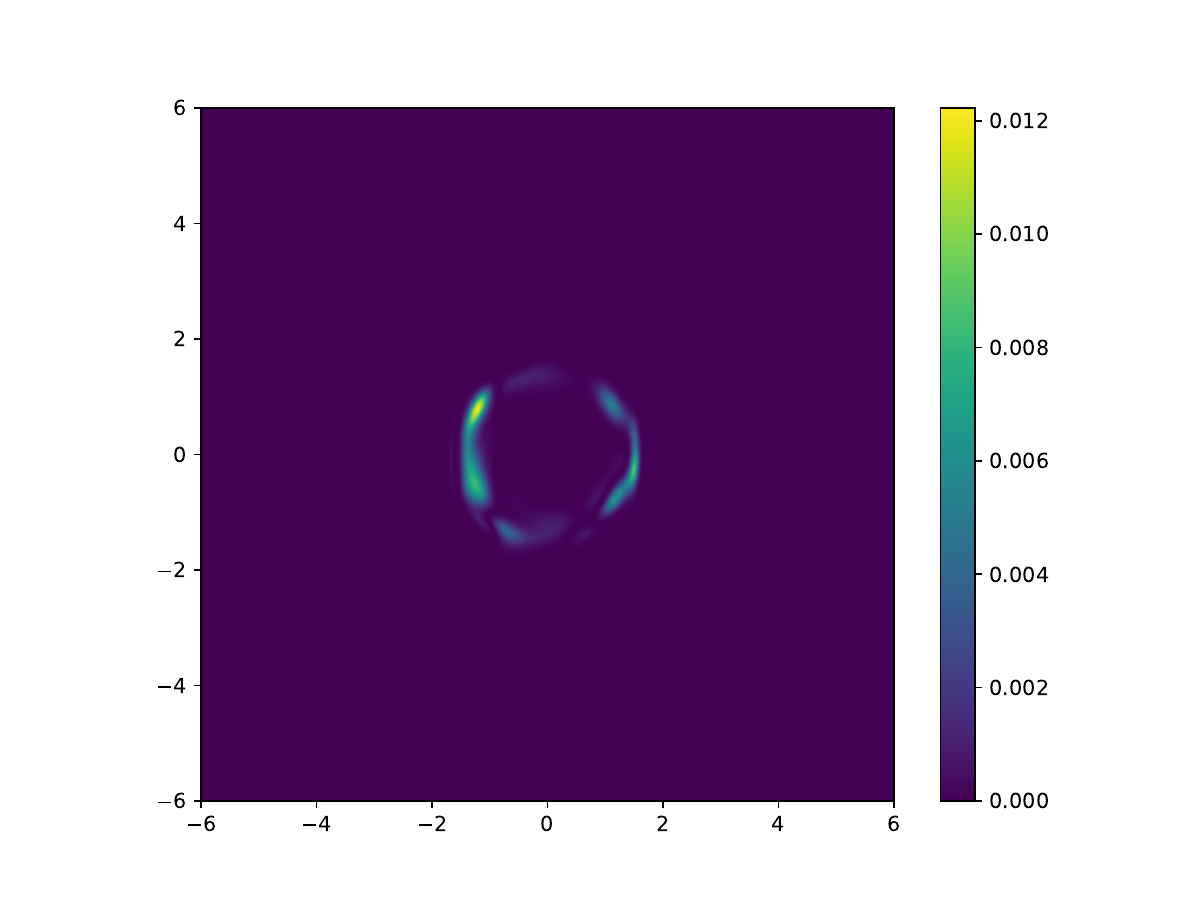}
    \caption{The squared residual norm evaluated after Bayesian training for the two-dimensional anisotropic, annular potential problem. The mapping illustrates regions where the constraints are less strictly satisfied, providing a visual proxy for the total predictive uncertainty across the domain.}
    \label{fig: 2D Anisotropic DoubleWell Resdiual - BNN}
\end{figure}

The neural network representation of the pullback map is obtained by training the model, which enables the computation of the densities by \eqref{eq:pullback-1}. The evolution of the data distribution and the reconstructed probability density functions using two distinct models are illustrated in Figure \ref{fig: 2D Anisotropic DoubleWell - PDF}. The probabilistic model in \eqref{eq:finprobloss} further enables uncertainty quantification. In particular, Bayesian neural networks provide estimates of the mean and variance of the reconstructed density, as shown in Figure \ref{fig: 2D Anisotropic DoubleWell - BNN}. As might be expected, the regions of high variance in this plot correspond to higher values of the residual component of the objective function \eqref{eq:optimization}, defined as $\left| \rho v + \rho \nabla \widetilde{\Psi} +\widetilde{D}\nabla \rho \right|^2$, as shown in Figure (\ref{fig: 2D Anisotropic DoubleWell Resdiual - BNN}).

\noindent\textbf{Remark}: This work utilizes a Bayesian framework, specifically Bayesian Neural Networks (BNNs), to provide robust  quantification of the uncertainty arising due to sparseness of the data. It is important to clarify that the observed residuals integrate uncertainty from three interconnected components: the transport map, the potential function, and the diffusion tensor. Due to the high degree of coupling between these variables, a rigorous causal attribution of uncertainty of each individual source remains a non-trivial challenge. We also note that variational inference, which has been employed in this work, allows a formal decomposition into the aleatoric and epistemic components of uncertainty  using approaches such as  Variational Uncertainty Decomposition \cite{jayasekera2025variational}. Such a deeper treatment of uncertainty is reserved as an important objective for future study to further enhance model interpretability.

\subsection{A high-dimensional Fokker-Planck problem}
\label{sec:highdim}
In this section, we apply our method to Fokker-Planck dynamics in five-dimensions. A primary motivating example is the analysis of experimental data on cell dynamics of molecular signaling and migration. Previous studies have demonstrated that cell migration can be modeled using variants of the Fokker-Planck equation \cite{srivastava2025inference}. Cell signaling involves distinct chemical reporters interacting with each other and the spatial position variables of migration. The number of these reporters can range from one to tens. With such an extension in mind, we demonstrate the convergence of the inferred diffusion tensor components and parameters of  the structured potential representation, as well as present projections of the inferred neural network potential surface. In the context of high-dimensional system inference, Chen et al. previously proposed an algorithm using the general Physics-Informed Neural Network framework, which captures both the governing equation and the underlying density function \cite{Chen2021a}, also in five dimensions. Their approach involves defining a loss function based on the PDE residual evaluated at collocation points. In similar numerical experiments, they used $100,000$ samples for each of seven time snapshots, imposed the PDE constraint using mini-batches of size $50,000$, and evaluated the residual on $10,000$ points. We show that our method achieves comparable convergence with significantly fewer data points.

\begin{figure}[ht!]
     \centering
     \includegraphics[width=0.6\textwidth]{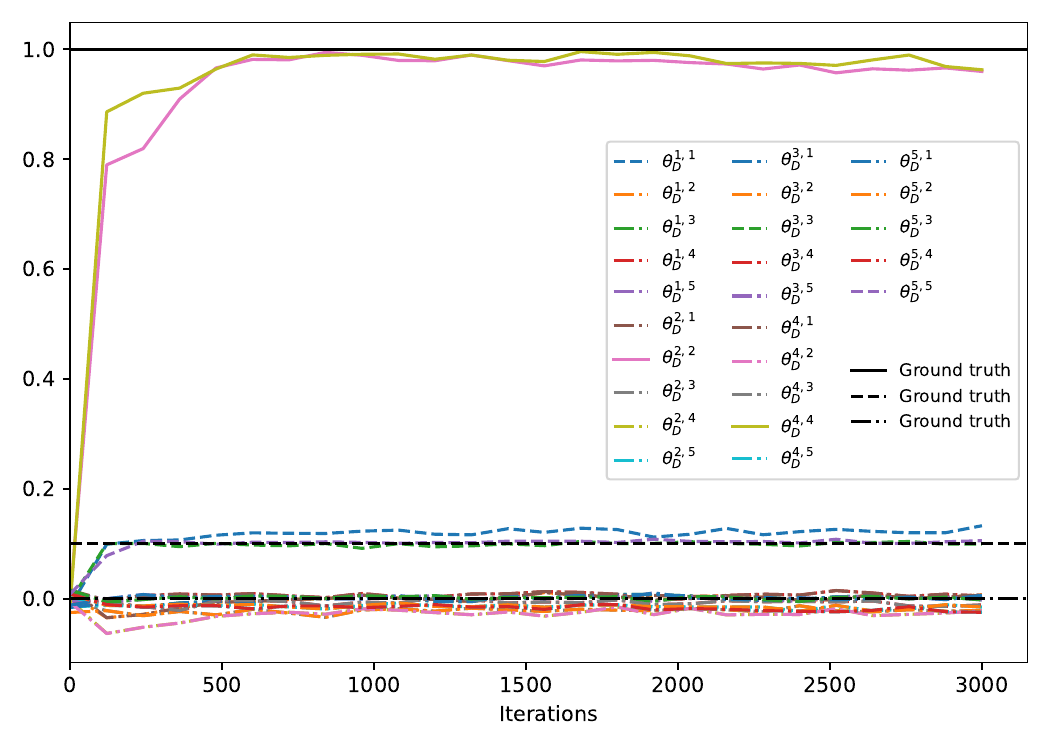}
     \includegraphics[width=0.6\textwidth]{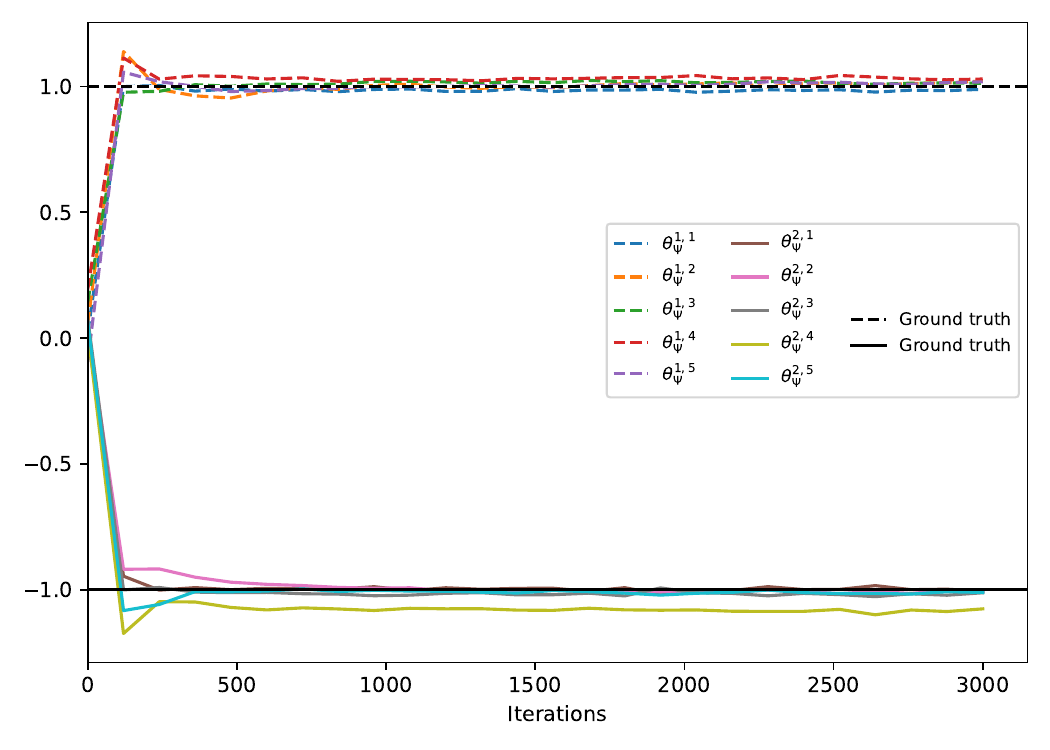}
    \caption{Convergence of the PDE system parameters in the five-dimensional anisotropic diffusion problem with a non-polynomial potential function. Top: diffusion tensor; bottom: potential function parameters over training iterations.}
    \label{fig: 5D Convergence}
\end{figure}

\begin{figure}[ht!]
     \centering
     \includegraphics[width=0.6\textwidth]{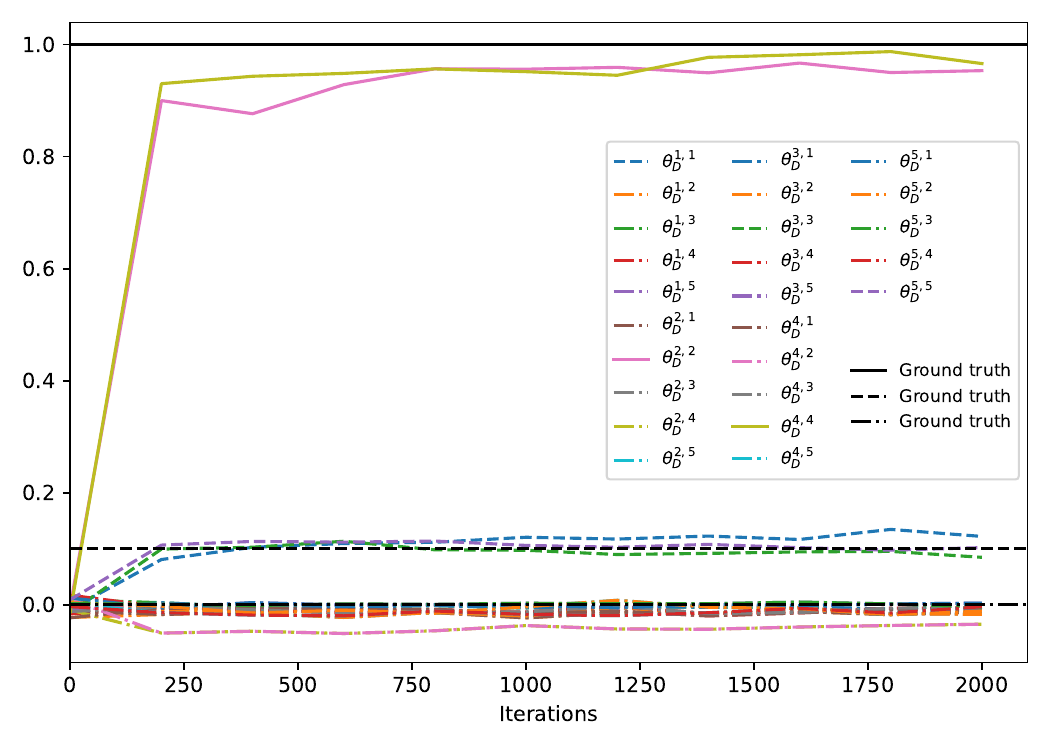}
    \caption{Convergence of the diffusion tensor parameters in the five-dimensional anisotropic diffusion problem with a non-polynomial potential function, represented using a neural network potential model.}
    \label{fig: 5D Convergence_2}
\end{figure}

\begin{figure}[ht!]
    \includegraphics[width=0.85\textwidth]{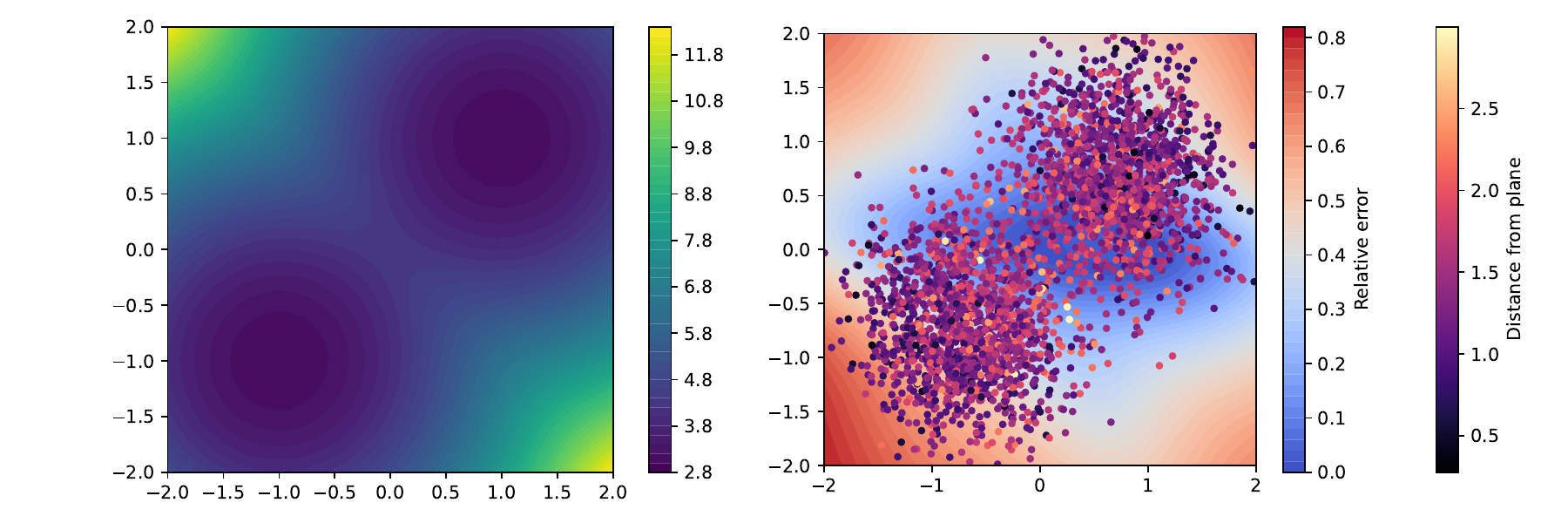}
    \includegraphics[width=0.85\textwidth]{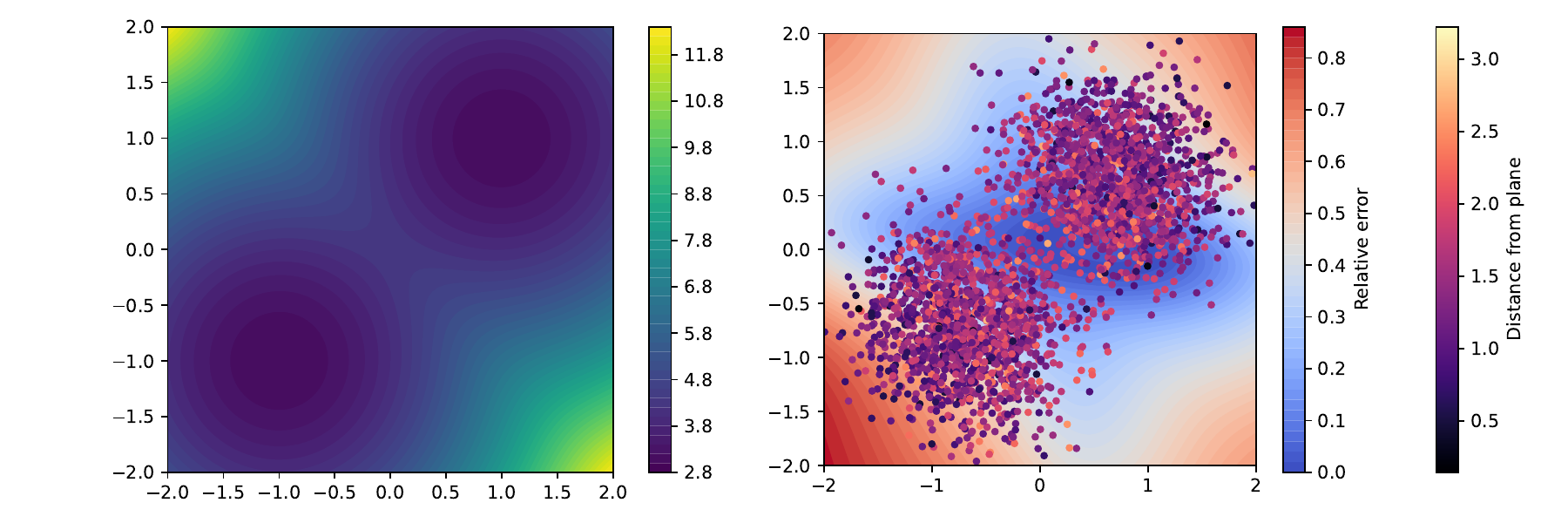}
    \includegraphics[width=0.85\textwidth]{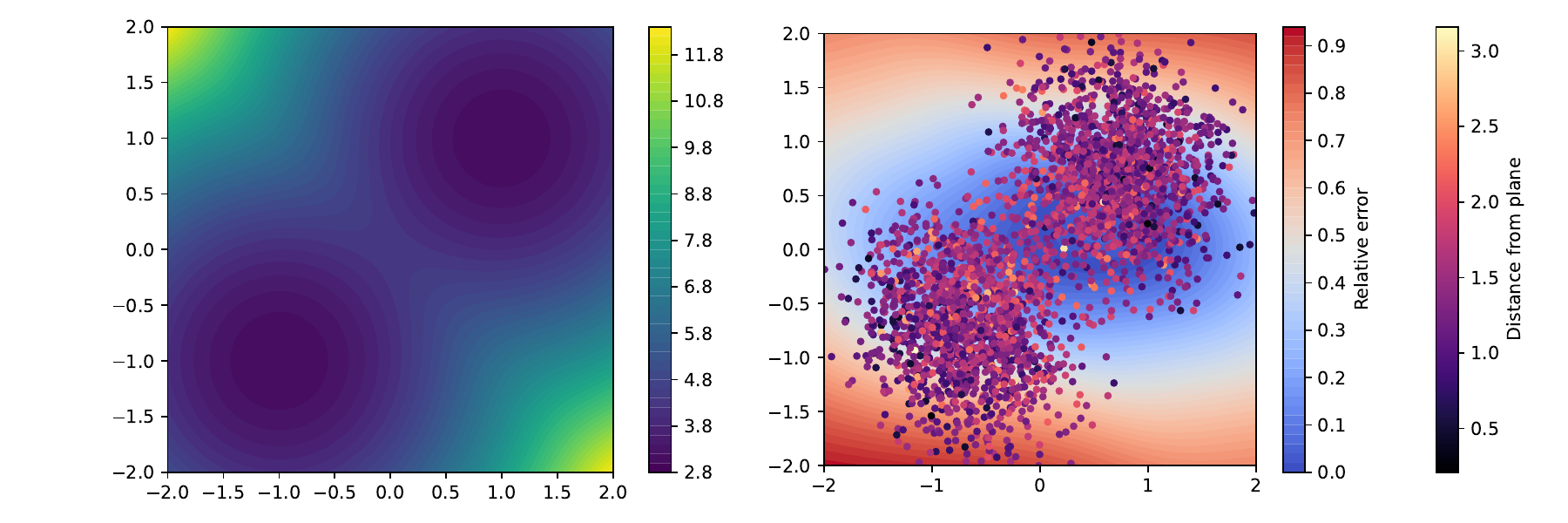}
    \caption{Each row displays a two-dimensional slice of the five-dimensional Fokker-Planck potential function, with the unplotted coordinates fixed at zero. From top to bottom: (i) $x_3=x_4=x_5=0$; (ii) $x_2=x_4=x_5=0$; (iii) $x_2=x_3=x_5=0$. Left: Contour plots of the ground truth Fokker-Planck potential function. Right: Contour plots of the relative error of the recovered neural-network potential function. Superposed on these error plots are data samples collected at time $t=0.5$ and projected onto the sliced plane. The color of each sample point indicates its distance from the respective two-dimensional slicing plane. }
    \label{fig: 5D Anisotropic NonPoly - Potential and Samples}
\end{figure}

\begin{figure}[ht!]
    \centering
    \includegraphics[width=0.85\textwidth]{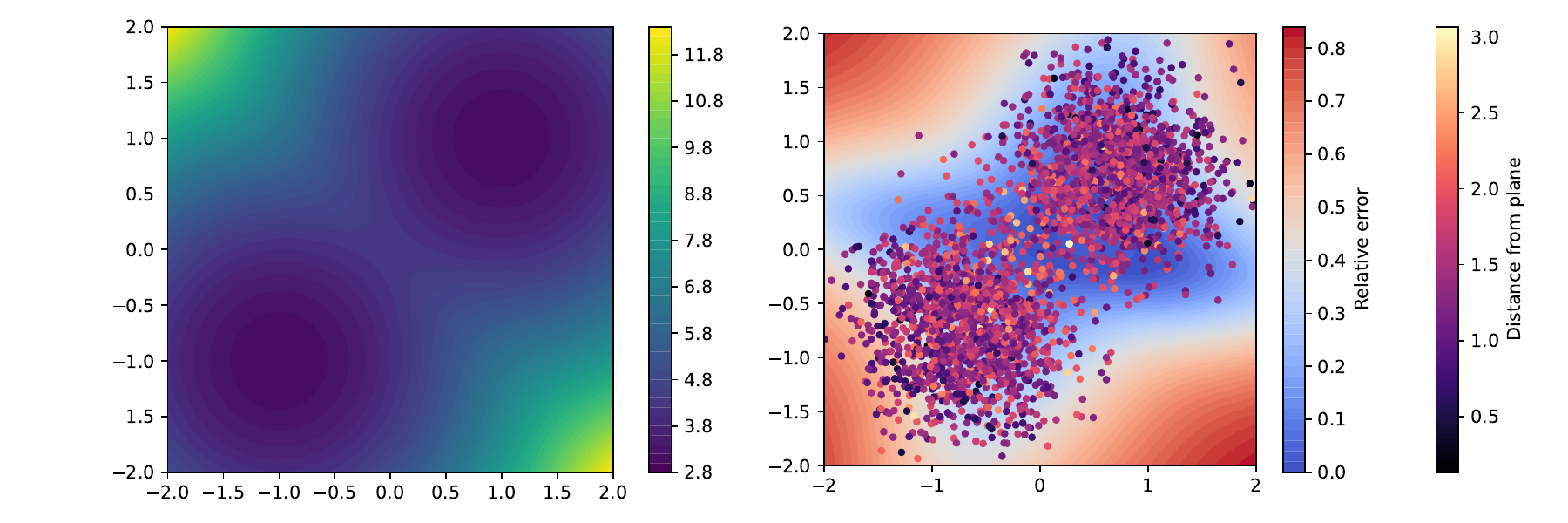}
    \includegraphics[width=0.85\textwidth]{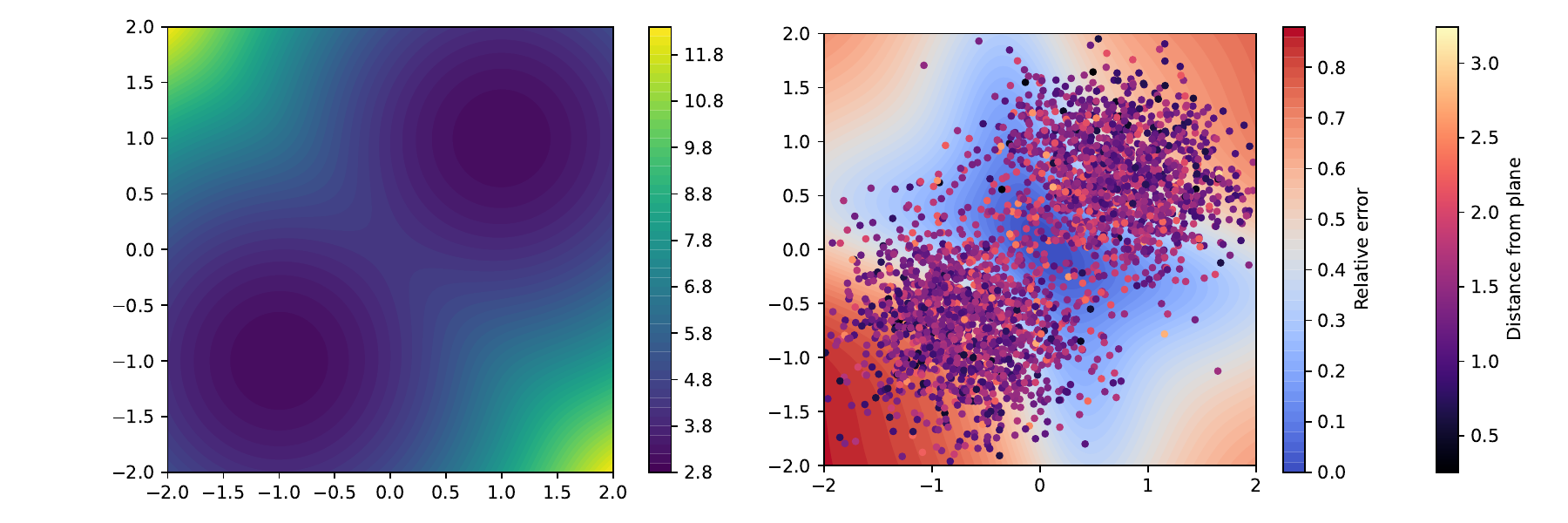}
    \includegraphics[width=0.85\textwidth]{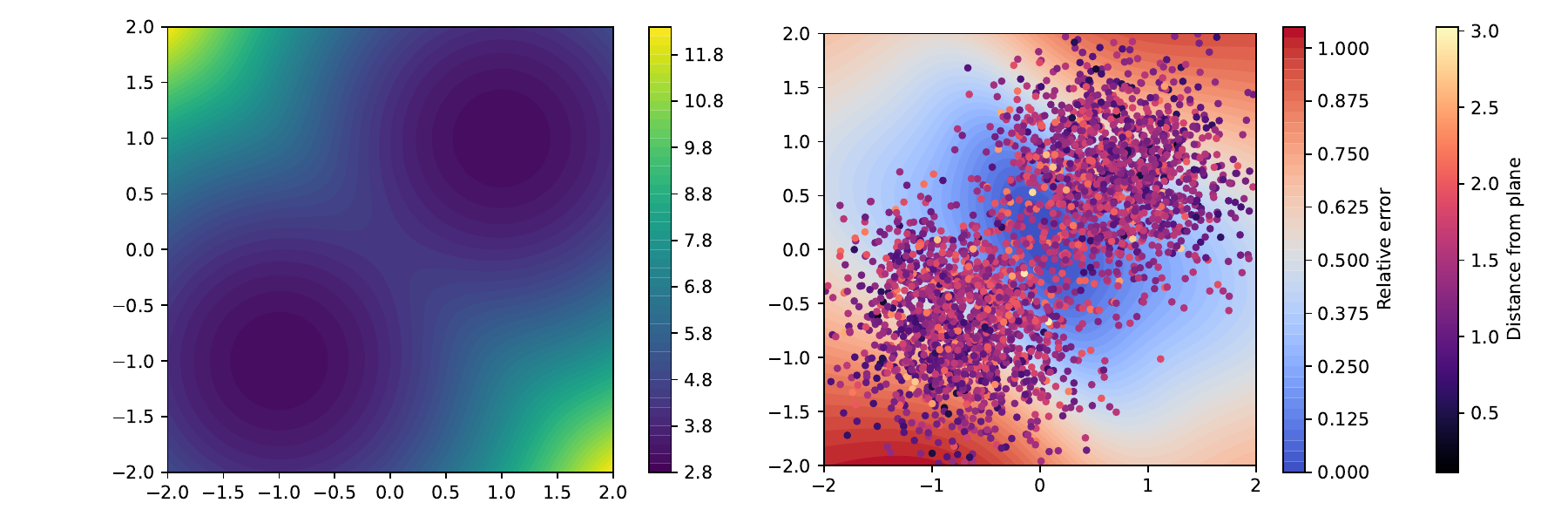}
    \caption{Each row displays a two-dimensional slice of the five-dimensional Fokker-Planck potential function, with the unplotted coordinates fixed at zero. From top to bottom: (i) $x_2=x_3=x_4=0$; (ii) $x_1=x_4=x_5=0$; (iii) $x_1=x_3=x_5=0$. Left: Contour plots of the ground truth Fokker-Planck potential function. Right: Contour plots of the relative error of the recovered neural-network potential function. Superposed on these error plots are data samples collected at time $t=0.5$ and projected onto the sliced plane. The color of each sample point indicates its distance from the respective two-dimensional slicing plane.}
    \label{fig: 5D Anisotropic NonPoly - Potential and Samples - 2}
\end{figure}

\begin{figure}[ht!]
\includegraphics[width=0.85\textwidth]{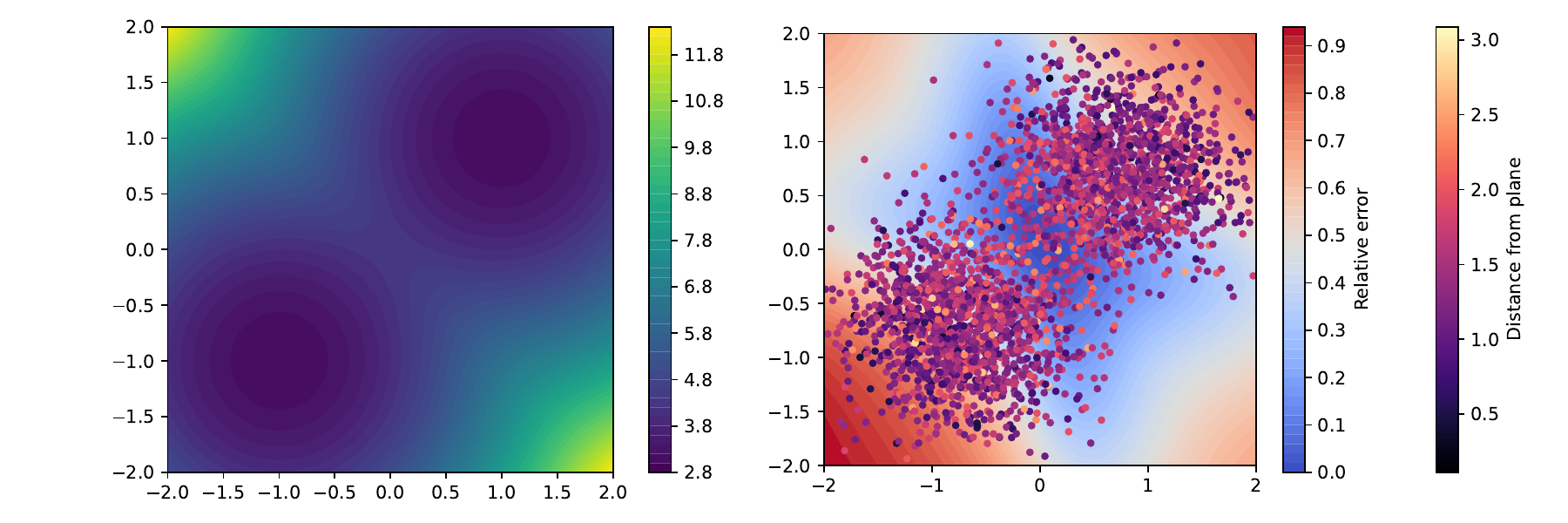}
\includegraphics[width=0.85\textwidth]{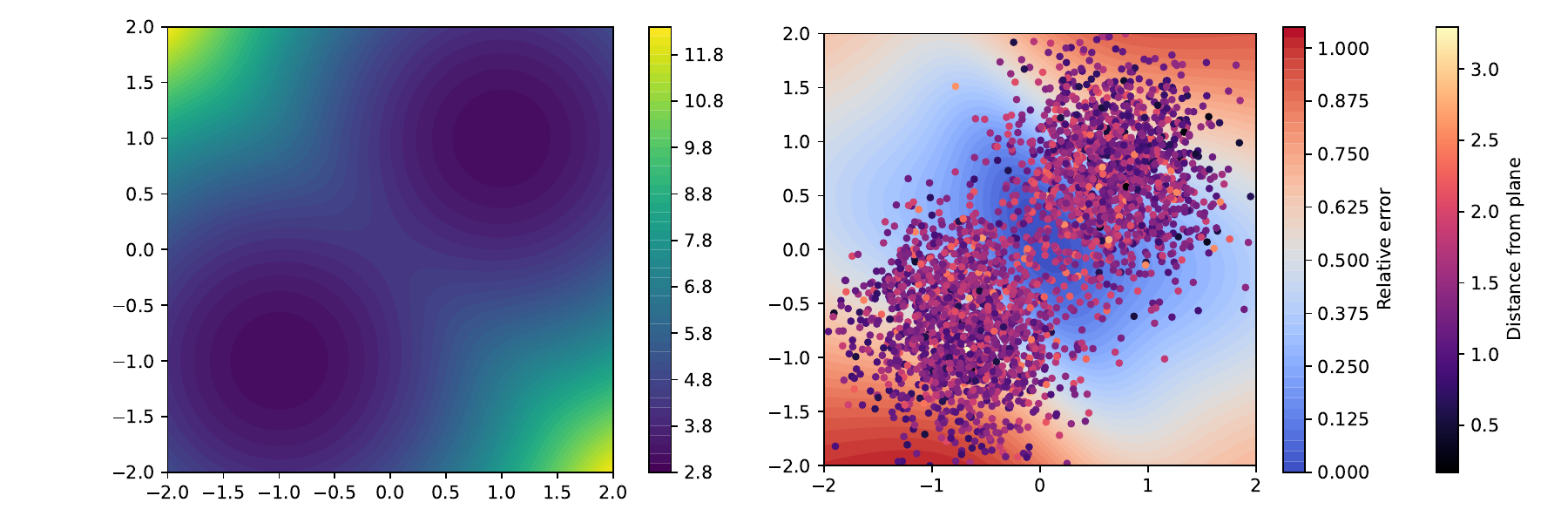}
\includegraphics[width=0.85\textwidth]{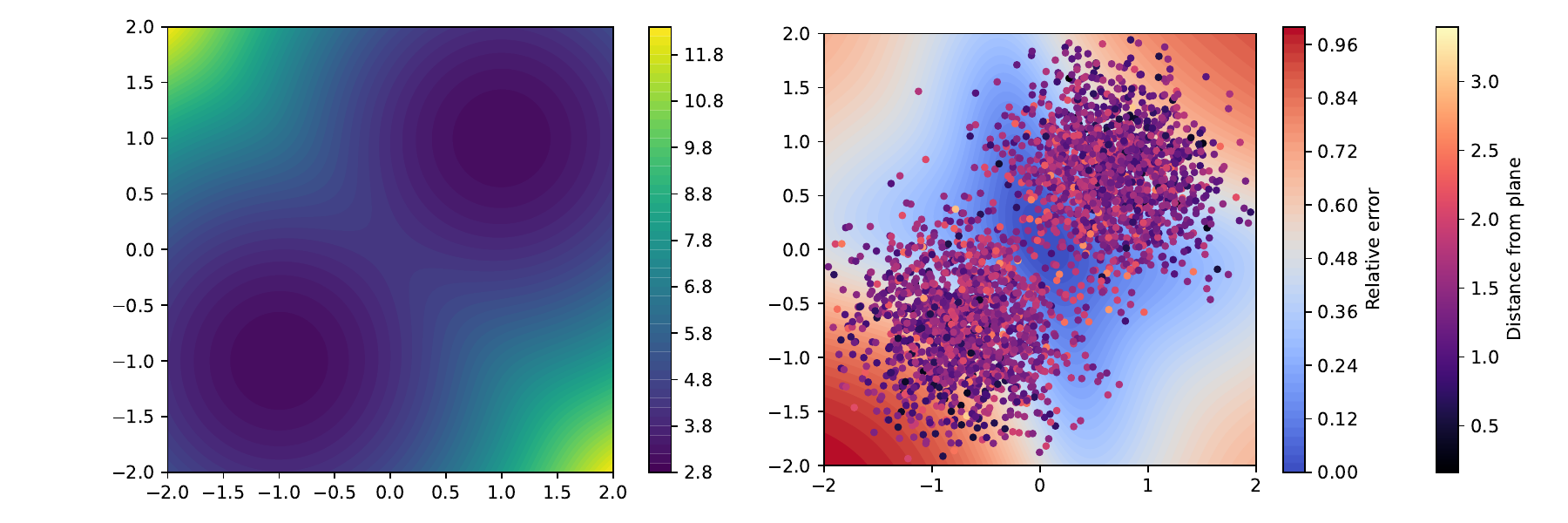}
\includegraphics[width=0.85\textwidth]{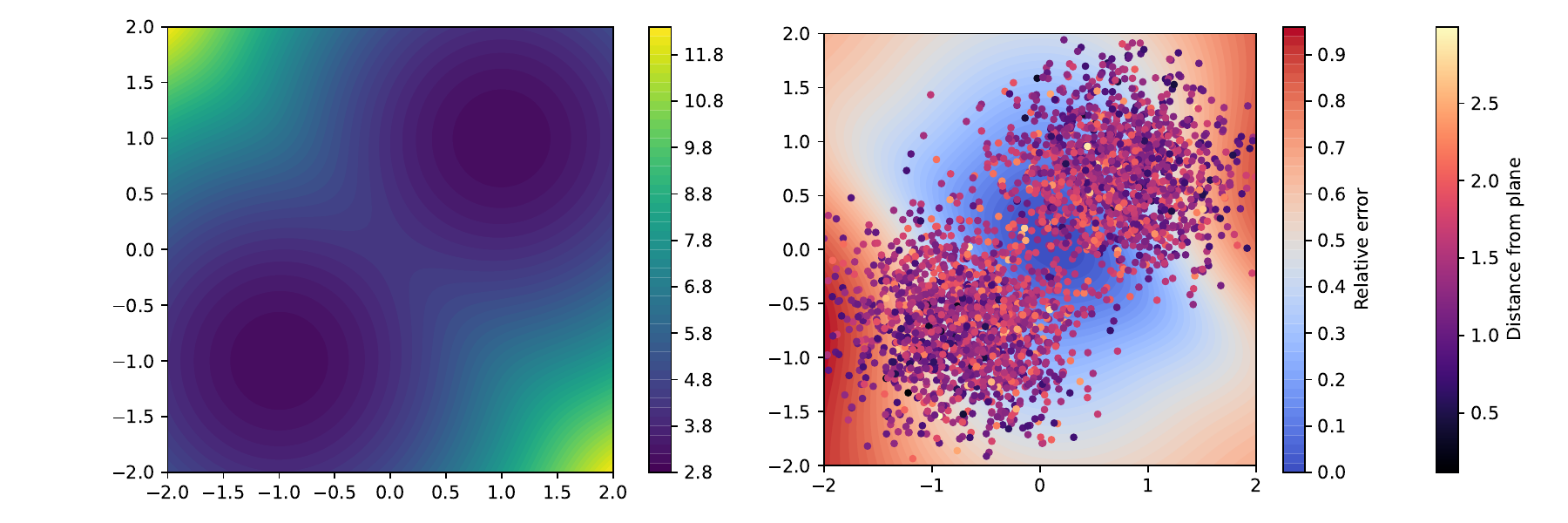}
 \caption{Each row displays a two-dimensional slice of the five-dimensional Fokker-Planck potential function, with the unplotted coordinates fixed at zero. From top to bottom: (i) $x_1=x_3=x_4=0$; (ii) $x_1=x_2=x_5=0$; (iii) $x_1=x_2=x_4=0$; (iv) $x_1=x_2=x_3=0$. Left: Contour plots of the ground truth Fokker-Planck potential function. Right: Contour plots of the relative error of the recovered neural-network potential function. Superposed on these error plots are data samples collected at time $t=0.5$ and projected onto the sliced plane. The color of each sample point indicates its distance from the respective two-dimensional slicing plane.}
    \label{fig: 5D Anisotropic NonPoly - Potential and Samples - 3}
\end{figure}

Consider the dynamics where the diffusion coefficients vary by a factor of $10$ across different directions, and the potential function is given in a non-polynomial form:
\begin{equation}
    D = \begin{pmatrix}
        0.1 & 0 & 0 & 0 & 0\\
        0 & 1 & 0 & 0 & 0\\
        0 & 0 & 0.1 & 0 & 0\\
        0 & 0 & 0 & 1 & 0 \\
        0 & 0 & 0 & 0 & 0.1
    \end{pmatrix}, \quad \Psi = -\log \left(e^{-\sum_{i=1}^5 (x_i + 1)^2} + e^{-\sum_{i=1}^5 (x_i - 1)^2}\right).
\end{equation}
We initialize the system with a standard Gaussian distribution and draw $3,000$ samples from it--almost two orders of magnitude fewer than used by Chen et al. \cite{Chen2021a}. To construct the dataset, we record the evolution of the system at time steps $t_i = 0.1i$, for $i=0,1,\cdots, 10$. The pullback maps are modeled using neural networks, and the PDE parameters are represented in the following form, noting that this is the structure-informed representation of the potential:  
\begin{equation}
    D = \begin{pmatrix}
        \theta_D^{1,1} & \theta_D^{1,2} & \theta_D^{1,3} & \theta_D^{1,4} & \theta_D^{1,5}\\
        \theta_D^{1,2} & \theta_D^{2,2} & \theta_D^{2,3} & \theta_D^{2,4} & \theta_D^{2,5}\\
        \theta_D^{1,3} & \theta_D^{2,3} & \theta_D^{3,3} & \theta_D^{3,4} & \theta_D^{3,5}\\
        \theta_D^{1,4} & \theta_D^{2,4} & \theta_D^{3,4} & \theta_D^{4,4} & \theta_D^{4,5}\\
        \theta_D^{1,5} & \theta_D^{2,5} & \theta_D^{3,5} & \theta_D^{4,5} & \theta_D^{5,5}
    \end{pmatrix}, \quad \Psi = -\log \left(e^{-\sum_{i=1}^5 (x_i + \theta_\Psi^{1,i})^2} + e^{-\sum_{i=1}^5 (x_i + \theta_\Psi^{2,i})^2}\right).
\end{equation}
Additionally, we used the neural network representation of the potential. Since the diffusion tensor must be symmetric, we enforce this by expressing it as $D = \frac{1}{2}(\tilde{D}+\tilde{D}^T)$. (Alternatively, one may parametrize only the lower triangular portion of the matrix to improve computational efficiency.) We track the convergence of the system parameters (diffusion tensor components and structured potential representation) over iterations, with the results shown in Figure \ref{fig: 5D Convergence} for the structure-informed potential function representation to allow comparison with the ground-truth synthetic data. Comparable results were obtained for the diffusion tensor when the full neural network representation was used for the potential with $12,000$ particles. Importantly, the transport map representation is always with neural networks. The results obtained using the neural network potential representation are illustrated in Figure \ref{fig: 5D Convergence_2} for convergence of the diffusion tensor, and  in Figures \ref{fig: 5D Anisotropic NonPoly - Potential and Samples}, \ref{fig: 5D Anisotropic NonPoly - Potential and Samples - 2}, \ref{fig: 5D Anisotropic NonPoly - Potential and Samples - 3} , which show two-dimensional slices through the surface that is the relative  error in the neural network potential. The ground-truth data on particle positions projected onto the corresponding plane at $t = 0.5$ is super-posed on the relative error in the neural network potential surface. The particle color represents the projection distance. While each two-dimensional slice of the ground-truth potential is identical by symmetry, this constraint was not imposed on the structured and neural network representations. The learnt neural network representations differ between two-dimensional slices, because the diffusional motion of $12,000$ particles over $\mathbb{R}^5$ provides sparse  information and unsymmetric particle distributions to resolve the surface. This results in different relative errors in the neural network potential surface for each slice. However, the minima within the slices are in good approximation. With lower diffusion, the samples would be more localized to the minima. 

\section{Discussion and conclusion}

In this paper, we have presented a framework for inferring  the   Fokker-Planck equation underlying the dynamics of an observed data distribution, while simultaneously learning an explicit transport map-based representation of the corresponding probability densities. This method is built upon the Knothe-Rosenblatt rearrangement, a triangular transport map that is both computationally tractable and theoretically well-established. In this regard, we note that  the variable ordering in the Rosenblatt transformation can influence the results. Generally, there is no universal criterion for determining an optimal coordinate sequence. The ideal permutation is often problem-dependent and requires domain-specific insights or additional structural assumptions. While developing systematic strategies for identifying optimal ordering is a compelling direction for future research, it remains beyond the scope of this study.

\begin{figure}[h]
     \centering
     \includegraphics[width=0.6\textwidth]{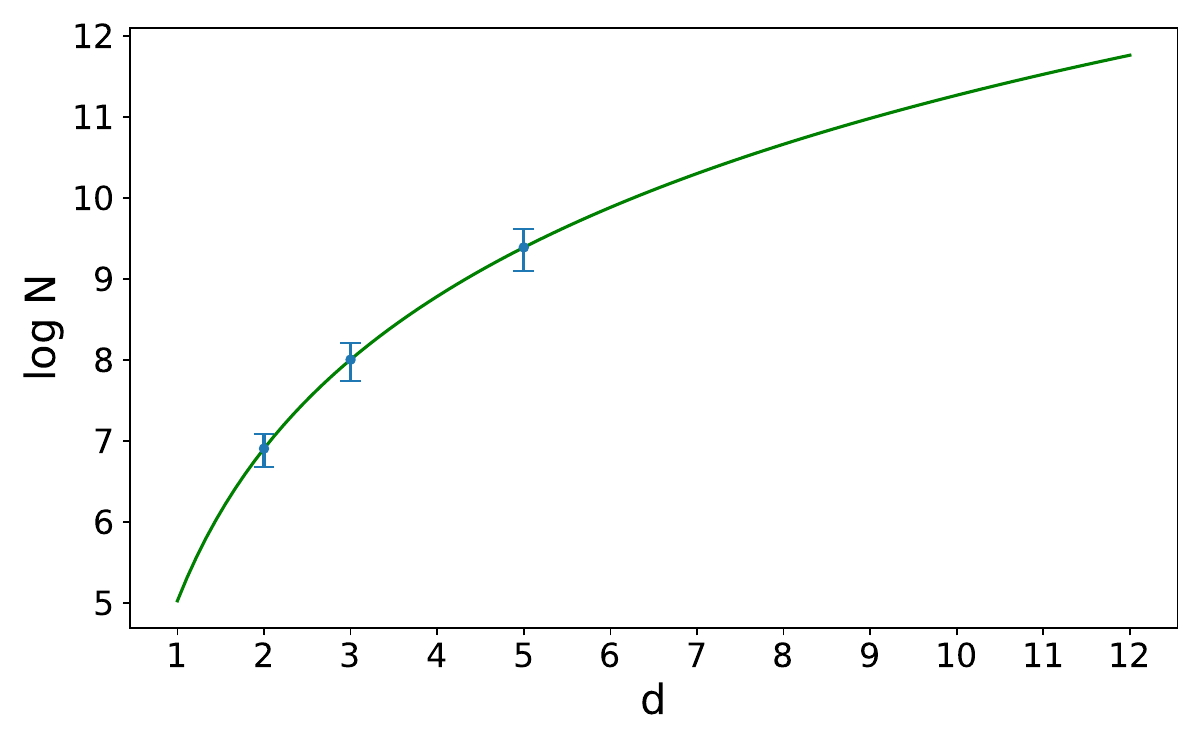}
    \caption{Empirical assessment of sample requirements. The plot illustrates the number of training samples required to achieve a stable reconstruction, demonstrating the model's scaling behavior. The empirical scaling is of the form $N = 152.93\times d^{2.71}$.}
    \label{fig: d_vs_logN}
\end{figure}

A key observation is that our approach requires far fewer data points to achieve robust performance, even compared to nontraditional grid-free methods. Empirically, in the absence of prior structural knowledge, approximately $1,000$ samples per time step were needed in two-dimensional settings to achieve a suitable accuracy with respect to ground truth data. Comparable results were obtained using around $3,000$, $7,000$, and $12,000$ samples per time step in three, four, and five dimensions, respectively, with motivation coming from cell dynamics in high dimensions. See Figure \ref{fig: d_vs_logN}, where we have included a fit of $\log\,N$ (samples) versus dimension $d$, which demonstrates $N = 152.93\times d^{2.71}$ scaling.
When prior structural information is incorporated, the required sample size decreased further; for instance only $3,000$ samples were needed in five dimensions. Even when integrals are estimated via random sampling, the transport-based representation adaptively assigns greater weight to regions where the data is concentrated, without requiring explicit knowledge of the boundaries of the data distribution. The above reported sample numbers are lower by orders of magnitude when compared with other methods in the literature. As noted in Section \ref{sec:highdim}, the transport map-based approach is more efficient by at least an order of magnitude in comparison with the work of Chen et al. \cite{Chen2021}.

We treat the uncertainty introduced by finite data samples in a standard approach via Bayesian neural networks trained using variational inference. The framework of variational inference allows a  decomposition into the aleatoric and epistemic components of uncertainty  using approaches such as  Variational Uncertainty Decomposition \cite{jayasekera2025variational}, which we reserve for future work.

Regarding the model architecture, we focus on the trajectory of transport maps by incorporating the time variable directly into a single neural network. This design enables the model to learn a time-continuous representation of the transport map, allowing for robust computation of time derivatives regardless of the time step size in the input data. The transport map furnishes the velocity in the target space, using which we obtain the total flux in the conservation form of the continuity equation. We then equate the transport map flux to the Fokker-Planck flux by minimizing an equivalent Lagrangian. This improves the numerical stability of the training by circumventing the imposition of the divergence operator in the loss. It also contributes to robust inference of the transport map and Fokker-Planck parameters with relatively few samples in two through five dimensions, which is as far as we have tested, guided by comparisons in the literature.   

\section*{Code}
Our code is available at \url{github.com/mechanoChem/transportMaps}

\section*{Acknowledgments}
This work was partly supported by the University of Michigan Rackham Graduate School (SH) and by a W.M. Keck Foundation grant (KG). 

\clearpage

\bibliographystyle{unsrt}
\bibliography{System_Identification_via_Transport_Maps}

\begin{thebibliography}{10}

\bibitem{Brunton2016}
S.L. Brunton, J.L. Proctor, and J.N. Kutz.
\newblock {Discovering governing equations from data by sparse identification of nonlinear dynamical systems}.
\newblock {\em Proc. Natl. Acad. Sci.}, 113(15):3932--3937, 2016.

\bibitem{Schaeffer2017}
H.~Schaeffer.
\newblock {Learning partial differential equations via data discovery and sparse optimization}.
\newblock {\em Proc. R. Soc. A Math. Phys. Eng. Sci.}, 473(2197), 2017.

\bibitem{Wang2019a}
Z.~Wang, X.~Huan, and K.~Garikipati.
\newblock {Variational system identification of the partial differential equations governing the physics of pattern-formation: Inference under varying fidelity and noise}.
\newblock {\em Comput. Methods Appl. Mech. Eng.}, 356:44--74, 2019.

\bibitem{wang2020system}
Z.~Wang, X.~Zhang, G.H. Teichert, M.~Carrasco-Teja, and K.~Garikipati.
\newblock System inference for the spatio-temporal evolution of infectious diseases: Michigan in the time of {COVID}-19.
\newblock {\em Computational Mechanics}, 66:1153--1176, 2020.

\bibitem{Wang2021}
Z.~Wang, X.~Huan, and K.~Garikipati.
\newblock {Variational system identification of the partial differential equations governing microstructure evolution in materials: Inference over sparse and spatially unrelated data}.
\newblock {\em Comput. Methods Appl. Mech. Eng.}, 377:113706, 2021.

\bibitem{wang2021system}
Z.~Wang, M.~Carrasco-Teja, X.~Zhang, G.H. Teichert, and K.~Garikipati.
\newblock System inference via field inversion for the spatio-temporal progression of infectious diseases: Studies of {COVID}-19 in {M}ichigan and {M}exico.
\newblock {\em Archives of Computational Methods in Engineering}, 28:4283--4295, 2021.

\bibitem{nikolov2022ogden}
D.P. Nikolov, S.~Srivastava, B.A. Abeid, U.M. Scheven, E.M. Arruda, K.~Garikipati, and J.B. Estrada.
\newblock Ogden material calibration via magnetic resonance cartography, parameter sensitivity and variational system identification.
\newblock {\em Philosophical Transactions of the Royal Society A}, 380(2234):20210324, 2022.

\bibitem{Messenger2022}
D.A. Messenger and D.M. Bortz.
\newblock {Learning mean-field equations from particle data using WSINDy}.
\newblock {\em Phys. D Nonlinear Phenom.}, 439:133406, 2022.

\bibitem{Raissi2019b}
M.~Raissi, P.~Perdikaris, and G.E. Karniadakis.
\newblock {Physics-informed neural networks: A deep learning framework for solving forward and inverse problems involving nonlinear partial differential equations}.
\newblock {\em J. Comput. Phys.}, 378:686--707, 2019.

\bibitem{Karniadakis2021}
G.E. Karniadakis, I.G. Kevrekidis, L.~Lu, P.~Perdikaris, S.~Wang, and L.~Yang.
\newblock {Physics-informed machine learning}.
\newblock {\em Nat. Rev. Phys.}, 3(6):422--440, 2021.

\bibitem{lu2019deeponet}
L.~Lu, P.~Jin, and G.E. Karniadakis.
\newblock Deeponet: Learning nonlinear operators for identifying differential equations based on the universal approximation theorem of operators.
\newblock {\em arXiv preprint arXiv:1910.03193}, 2019.

\bibitem{lu2021learning}
L.~Lu, P.~Jin, G.~Pang, Z.~Zhang, and G.E. Karniadakis.
\newblock Learning nonlinear operators via deeponet based on the universal approximation theorem of operators.
\newblock {\em Nature machine intelligence}, 3(3):218--229, 2021.

\bibitem{li2020fourier}
Z.~Li, N.~Kovachki, K.~Azizzadenesheli, B.~Liu, K.~Bhattacharya, A.~Stuart, and A.~Anandkumar.
\newblock Fourier neural operator for parametric partial differential equations.
\newblock {\em arXiv preprint arXiv:2010.08895}, 2020.

\bibitem{li2023fourier}
Z.~Li, D.Z. Huang, B.~Liu, and A.~Anandkumar.
\newblock Fourier neural operator with learned deformations for pdes on general geometries.
\newblock {\em Journal of Machine Learning Research}, 24(388):1--26, 2023.

\bibitem{Garrison2008}
J.C. Garrison and R.Y. Chiao.
\newblock {The master equation}.
\newblock {\em Quantum Opt.}, pages 538--577, 2008.

\bibitem{risken1989fokker}
H.~Risken.
\newblock Fokker-planck equation.
\newblock In {\em The Fokker-Planck equation: methods of solution and applications}, pages 63--95. Springer, 1989.

\bibitem{Erban2009}
R.~Erban and S.J. Chapman.
\newblock {Stochastic modelling of reaction-diffusion processes: Algorithms for bimolecular reactions}.
\newblock {\em Phys. Biol.}, 6(4), 2009.

\bibitem{Walczak2012}
A.M. Walczak, A.~Mugler, and C.H. Wiggins.
\newblock {Analytic methods for modeling stochastic regulatory networks}.
\newblock {\em Methods Mol. Biol.}, 880:273--322, 2012.

\bibitem{Sato2014}
I.~Sato and H.~Nakagavva.
\newblock {Approximation analysis of stochastic gradient langevin dynamics by using fokker-planck equation and ito process}.
\newblock {\em 31st Int. Conf. Mach. Learn. ICML 2014}, 3:2647--2655, 2014.

\bibitem{Dai2020}
X.~Dai and Y.~Zhu.
\newblock {\em On Large Batch Training and Sharp Minima: A Fokker–Planck Perspective}, volume~14.
\newblock Springer International Publishing, 2020.

\bibitem{jordan1998variational}
R.~Jordan, D.~Kinderlehrer, and F.~Otto.
\newblock The variational formulation of the fokker-planck equation.
\newblock {\em SIAM journal on mathematical analysis}, 29(1):1--17, 1998.

\bibitem{Liu2021}
W.~Liu, C.K.L. Kou, K.H. Park, and H.K. Lee.
\newblock {Solving the inverse problem of time independent Fokker–Planck equation with a self supervised neural network method}.
\newblock {\em Sci. Rep.}, 11(1):1--11, 2021.

\bibitem{Villani2003}
C.~Villani.
\newblock {\em {Topics In Optimal Transportation}}.
\newblock 2003.

\bibitem{Villani2008}
C.~Villani.
\newblock {Optimal Transport: Old and New}.
\newblock 2008.

\bibitem{Korotin2021a}
A.~Korotin, V.~Egiazarian, A.~Asadulaev, A.~Safin, and E.~Burnaev.
\newblock {Wasserstein-2 Generative Networks}.
\newblock {\em ICLR 2021 - 9th Int. Conf. Learn. Represent.}, (2016):1--30, 2021.

\bibitem{Gangbo1996}
W.~Gangbo and R.J. McCann.
\newblock {The geometry of optimal transportation}.
\newblock {\em Acta Math.}, 177(2):113--161, 1996.

\bibitem{Makkuva2020}
A.V. Makkuva, A.~Taghvaei, J.D. Lee, and S.~Oh.
\newblock {Optimal transport mapping via input convex neural networks}.
\newblock {\em 37th Int. Conf. Mach. Learn. ICML 2020}, PartF16814:6628--6637, 2020.

\bibitem{Knothe1957}
H.~Knothe.
\newblock {Contributions to the theory of convex bodies.}, 1957.

\bibitem{Rosenblatt1952}
M.~Rosenblatt.
\newblock {Remarks on a Multivariate Transformation}.
\newblock {\em Ann. Math. Stat.}, 23(3):470--472, 1952.

\bibitem{Santambrogio2015c}
F.~Santambrogio.
\newblock {\em {Optimal transport for applied mathematicians : calculus of variations, PDEs, and modeling}}.
\newblock Number May. 2015.

\bibitem{Kobyzev2021}
I.~Kobyzev, S.J.D Prince, and M.A. Brubaker.
\newblock {Normalizing Flows: An Introduction and Review of Current Methods}, 2021.

\bibitem{Tang2020}
Keju Tang, Xiaoliang Wan, and Qifeng Liao.
\newblock {Deep density estimation via invertible block-triangular mapping}.
\newblock {\em Theor. Appl. Mech. Lett.}, 10(3):143--148, 2020.

\bibitem{Baptista2023}
R.~Baptista, Y.~Marzouk, and O.~Zahm.
\newblock {On the Representation and Learning of Monotone Triangular Transport Maps}.
\newblock {\em Found. Comput. Math.}, 2023.

\bibitem{Li2017}
C.L. Li, W.C. Chang, Y.~Cheng, Y.~Yang, and B.~P{\'{o}}czos.
\newblock {MMD GAN: Towards deeper understanding of moment matching network}.
\newblock {\em Adv. Neural Inf. Process. Syst.}, 2017-Decem(MMD):2204--2214, 2017.

\bibitem{Arjovsky2017}
M.~Arjovsky, S.~Chintala, and L.~Bottou.
\newblock {Wasserstein GAN}.
\newblock 2017.

\bibitem{Yadav2018}
A.~Yadav, S.~Shah, Z.~Xu, D.~Jacobs, and T.~Goldstein.
\newblock {Stabilizing adversarial nets with prediction methods}.
\newblock {\em 6th Int. Conf. Learn. Represent. ICLR 2018 - Conf. Track Proc.}, 2018.

\bibitem{Seguy2018}
V.~Seguy and M.~Apping.
\newblock Large-scale optimal transport and mapping esetimation.
\newblock (1781):1--15, 2018.

\bibitem{Amos2017}
B.~Amos, L.~Xu, and J.Z. Kolter.
\newblock {Input convex neural networks: Supplementary material}.
\newblock {\em 34th Int. Conf. Mach. Learn. ICML 2017}, 1:192--206, 2017.

\bibitem{Brenier1991}
Y.~Brenier.
\newblock {Polar factorization and monotone rearrangement of vector‐valued functions}.
\newblock {\em Commun. Pure Appl. Math.}, 44(4):375--417, 1991.

\bibitem{Taghvaei2019}
A.~Taghvaei and A.~Jalali.
\newblock {2-Wasserstein Approximation via Restricted Convex Potentials with Application to Improved Training for GANs}.
\newblock pages 1--38, 2019.

\bibitem{McCann1995}
R.J. McCann.
\newblock {Existence and uniqueness of monotone measure-preserving maps}.
\newblock {\em Duke Math. J.}, 80(2):309--323, 1995.

\bibitem{Carlier2009}
G.~Carlier, A.~Galichon, and F.~Santambrogio.
\newblock {From Knothe's transport to Brenier's map and a continuation method for optimal transport}.
\newblock {\em SIAM J. Math. Anal.}, 41(6):2554--2576, 2009.

\bibitem{Blei2017}
D.M. Blei, A.~Kucukelbir, and J.D. McAuliffe.
\newblock {Variational Inference: A Review for Statisticians}.
\newblock {\em J. Am. Stat. Assoc.}, 112(518):859--877, 2017.

\bibitem{Graves2011}
A.~Graves.
\newblock {Practical variational inference for neural networks}.
\newblock {\em Adv. Neural Inf. Process. Syst. 24 25th Annu. Conf. Neural Inf. Process. Syst. 2011, NIPS 2011}, pages 1--9, 2011.

\bibitem{Zhang2023}
X.~Zhang and K.~Garikipati.
\newblock {Label-free learning of elliptic partial differential equation solvers with generalizability across boundary value problems}.
\newblock {\em Comput. Methods Appl. Mech. Eng.}, 417, 2023.

\bibitem{Talay1994}
D.~Talay.
\newblock {\em {Numerical solution of stochastic differential equations}}, volume~47.
\newblock 1994.

\bibitem{jayasekera2025variational}
I~Shavindra Jayasekera, Jacob Si, Filippo Valdettaro, Wenlong Chen, A~Aldo Faisal, and Yingzhen Li.
\newblock Variational uncertainty decomposition for in-context learning.
\newblock {\em arXiv preprint arXiv:2509.02327}, 2025.

\bibitem{srivastava2025inference}
Siddhartha Srivastava, Patrick~C Kinnunen, Zhenlin Wang, Kenneth~KY Ho, Brock~A Humphries, Siyi Chen, Jennifer~J Linderman, Gary~D Luker, Kathryn~E Luker, and Krishna Garikipati.
\newblock Inference of weak-form partial differential equations describing migration and proliferation mechanisms in wound healing experiments on cancer cells.
\newblock {\em PLOS Computational Biology}, 21(10):e1013607, 2025.

\bibitem{Chen2021a}
X.~Chen, L.~Yang, J.~Duan, and G.E. Karniadakis.
\newblock {Solving inverse stochastic problems from discrete particle observations using the fokker-planck equation and physics-informed neural networks}.
\newblock {\em SIAM J. Sci. Comput.}, 43(3):B811--B830, 2021.

\bibitem{Chen2021}
Zhao Chen, Yang Liu, and Hao Sun.
\newblock {Physics-informed learning of governing equations from scarce data}.
\newblock {\em Nat. Commun.}, 12(1):1--13, 2021.

\end{thebibliography}

\end{document}